\documentclass[final]{siamltex}

\usepackage{epsfig}
\title{A Refined Harmonic Lanczos Bidiagonalization Method and an
Implicitly Restarted Algorithm for Computing the Smallest Singular
Triplets of Large Matrices}

\author{Zhongxiao Jia\thanks{Department of Mathematical Sciences,
Tsinghua University, Beijing 100084, P. R. China (
jiazx@tsinghua.edu.cn,
http://faculty.math.tsinghua.edu/$^{\sim}$zjia). This author's
research was Supported by the NSFC Grant 10771116 and the Doctoral
Program of Ministry of Education (Grant 20060003003).} \and Datian
Niu\thanks{Department of Applied Mathematics, Dalian University of
Technology, Dalian 116024, P. R. China (niudt@dlnu.edu.cn). This
author's research was Supported by the NSFC Grants 10872045,
60501021.}}

\begin{document}
\maketitle
\begin{abstract}
The harmonic Lanczos bidiagonalization method can be used to compute
the smallest singular triplets of a large matrix $A$. We prove that
for good enough projection subspaces harmonic Ritz values converge
if the columns of $A$ are strongly linearly independent. On the
other hand, harmonic Ritz values may miss some desired singular
values when the columns of $A$ almost linearly dependent.
Furthermore, harmonic Ritz vectors may converge irregularly and even
may fail to converge. Based on the refined projection principle for
large matrix eigenproblems due to the first author, we propose a
refined harmonic Lanczos bidiagonalization method that takes the
Rayleigh quotients of the harmonic Ritz vectors as approximate
singular values and extracts the best approximate singular vectors,
called the refined harmonic Ritz approximations, from the given
subspaces in the sense of residual minimizations. The refined
approximations are shown to converge to the desired singular vectors
once the subspaces are sufficiently good and the Rayleigh quotients
converge. An implicitly restarted refined harmonic Lanczos
bidiagonalization algorithm (IRRHLB) is developed. We study how to
select the best possible shifts, and suggest refined harmonic shifts
that are theoretically better than the harmonic shifts used within
the implicitly restarted Lanczos bidiagonalization algorithm
(IRHLB). We propose a novel procedure that can numerically compute
the refined harmonic shifts efficiently and accurately. Numerical
experiments are reported that compare IRRHLB with five other
algorithms based on the Lanczos bidiagonalization process. It
appears that IRRHLB is at least competitive with them and can be
considerably more efficient when computing the smallest singular
triplets.

\end{abstract}



\pagestyle{myheadings} \thispagestyle{plain}

\markboth{ZHONGXIAO JIA AND DATIAN NIU}{A REFINED HARMONIC LANCZOS
BIDIAGONALIZATION METHOD}

\section{Introduction}

We assume that a large sparse matrix $A \in {\cal R}^{M \times
N},\,M\geq N$ has full column rank and let
\begin{equation}\label{svd}
    A=U\left ( \Sigma \atop 0 \right )V^{\rm T}=U_1 \Sigma V^{\rm T}
\end{equation}
be its singular value decomposition (SVD) \cite{golub96,stewart01},
where $U=(u_1,u_2,\ldots,u_M)=(U_1,U_2)$ and
$V=(v_1,v_2,\ldots,v_N)$ are $M\times M$ and $N\times N$ orthogonal
matrices, $U_1=(u_1,u_2,\ldots,u_N)$ and
$\Sigma=diag(\sigma_1,\sigma_2,\ldots,\sigma_N)$ is diagonal.
${\sigma}_i,i=1,2,\ldots,N$, are called the singular values of $A$,
$u_i$'s and $v_i$'s are the associated left and right singular
vectors, respectively, and $(\sigma_i,u_i,v_i)$'s are called
singular triplets. In this paper, slightly different from the
convention, the singular values are labeled as
$\sigma_1\leq\sigma_2\leq\cdots\leq\sigma_N$.

We are concerned with the following problem.

{\bf Problem 1}. {\em Compute numerically the $k$ smallest
singular triplets $(\sigma_i,u_i,v_i)$ of $A$, $i=1,2,\ldots,k$,
where $k\ll N$.}

There are many applications of Problem 1, including determination of
numerical rank and of spectral condition number, least squares
problems, total least squares problems, regression analysis, image
and signal processing, pattern recognition and information
retrieval, to name a few.

Consider the $(M+N)\times (M+N)$ augmented matrix
\begin{equation}\label{argumented}
    \tilde A=\left ( \begin{array} {c c}
                           0   & A \\
                           A^{\rm T} & 0
                     \end{array} \right ).
\end{equation}
Then, the eigenvalues of $\tilde A$ are just $\pm
\sigma_1,\ldots,\pm \sigma_N$ and $M-N$ zeros, the associated
eigenvectors of $\sigma_i$ and $-\sigma_i$ are
$\frac{1}{\sqrt{2}}\left (u_i^{\rm T},v_i^{\rm T} \right )^{\rm T}$
and $\frac{1}{\sqrt{2}} \left (u_i^{\rm T},-v_i^{\rm T} \right
)^{\rm T}$, respectively, and the eigenvectors associated with zero
eigenvalues have the form $\left(u^{\rm T},0^{\rm T} \right )^{\rm
T}$, where $u$'s are orthogonal to all $u_1,\ldots,u_N$. Therefore,
we obtain the following formulation of Problem 1.

{\bf Problem 2}. {\em Compute numerically the $k$ smallest
positive eigenvalues and the associated eigenvectors of $\tilde
A$.}

For the $k$ smallest eigenpairs of $\tilde A$, Problem 2 is a
symmetric interior eigenvalue problem. Since $M$ and $N$ are assumed
to be large, we can only resort to projection methods. A typical
method is the symmetric Lanczos method \cite{parlett}. It and other
standard projection methods usually favor the extreme eigenvalues
and the associated eigenvectors but are generally very inefficient
for computing interior eigenpairs \cite{parlett}. Another drawback
is that in finite precision the computed eigenvalues do not come in
plus-and-minus pairs and the computed eigenvectors do not respect
the special structures that the true eigenvectors have.

Because of the mentioned drawbacks, we should not work on $\tilde A$
explicitly for computing the smallest singular triplets of $A$.
Instead we attempt to solve Problem 1 directly by working on $\tilde
A$ implicitly. It appears that Lanczos bidiagonalization type
methods \cite{Baglama05,Hochstenbach2001,
Hochstenbach2002,jianiu03,Kokio04,larsen1,Simon00} and
Jacobi-Davidson SVD type methods \cite{Hochstenbach2001,
Hochstenbach2002} can solve the mentioned problems elegantly. The
Lanczos bidiagonalization type methods available have in common that
they are all based on the Lanczos bidiagonalization process to build
up orthonormal bases of certain Krylov subspaces. However, their
mathematical backgrounds can be fundamentally different. Basically,
there are three kinds of projection principles that extract
different approximate singular triplets with respect to the
subspaces. Some methods use the standard projection principle
\cite{book,parlett,stewart01} to extract Ritz approximations
\cite{Baglama05,baglama06,hernandez,Hochstenbach2001,Hochstenbach2002,jianiu03,
larsen1,Simon00,stoll}, some methods use the harmonic projection
principle \cite{book,stewart01,vorst} to extract harmonic Ritz
approximations
\cite{Baglama05,baglama06,Hochstenbach2001,Hochstenbach2002,Kokio04}
and some methods use the refined projection principle
\cite{book,jia97,stewart01,vorst} to extract refined singular vector
approximations \cite{Hochstenbach2002,jianiu03,Kokio04}.
Jacobi-Davidson type SVD methods for Problem 1 have several versions
that are based on the three projection principles as well as their
generalizations, respectively. As observed and claimed in
\cite{Baglama05,Hochstenbach2002}, the refined extraction version
appears to give the best accuracy in general.

For Problem 1, due to the storage requirement and computational
cost, all the Lanczos bidiagonalization type methods as well as
Jacobi-Davidson type methods have to be restarted generally  in
order to make them converge. That is, for given projection
subspaces, if the methods do not converge, then one repeatedly
chooses new better starting vectors, constructs better subspaces and
computes new approximate singular triplets until they converge. The
implicit restarting technique due to Sorensen \cite{sorensen92} is a
powerful tool for restarting Krylov subspace algorithms in various
contexts including large SVD problems
\cite{Baglama05,baglama06,Hochstenbach2001,
Hochstenbach2002,jianiu03,Kokio04,larsen1,Simon00}. The success of
an implicitly restarted algorithm heavily depends on both the
underlying method itself and a proper selection of the shifts
involved; see, e.g., \cite{jia99a,sorensen92}. Based on the Lanczos
bidiagonalization method and one of its harmonic versions, Jia and
Niu \cite{jianiu03} and Larsen \cite{larsen1} have developed an
implicitly restarted Lanczos bidiagonalization algorithm (IRLB), and
Kokiopoulou {\em et al.} \cite{Kokio04} have proposed an implicitly
restarted harmonic Lanczos bidiagonalization algorithm (IRLANB) for
computing the smallest singular triplets. IRLB uses the unwanted
Ritz values and IRLANB uses the unwanted Ritz or harmonic Ritz
values as shifts, respectively. These shifts are called exact shifts
and harmonic shifts. Baglama and Reichel \cite{Baglama05,baglama06}
propose a thick restarting technique that explicitly augments small
subspaces with certain Ritz or harmonic Ritz vectors, leading to
augmented restarted Lanczos bidiagonalization algorithms (IRLBA).
Hernandez {\em et al.} \cite{hernandez} analyze a parallel
implementation of this algorithm. Based on Stewart's work for large
eigenproblems \cite{stewart}, Stoll \cite{stoll} presents a
Krylov-Schur type algorithm that is restarted explicitly and is
easily implemented.

It is shown in \cite{jianiu03} that the Lanczos bidiagonalization
method may fail to compute singular vectors though it converges
for computing singular values for sufficiently good subspaces. To
correct this deficiency, applying the refined projection principle
proposed by the first author \cite{jia97} (see also
\cite{book,stewart01,vorst}), we have proposed a refined Lanczos
bidiagonalization method, analyzed its convergence and developed
an implicitly restarted refined Lanczos bidiagonalization
algorithm (IRRLB) \cite{jianiu03}. Based on the refined
approximations to singular vectors, we have proposed refined
shifts that are theoretically better than the exact shifts used
within IRLB. Numerical experiments have demonstrated that IRRLB
often outperforms IRLB \cite{jianiu03,larsen1} considerably and is
more efficient than several other available schemes: PROPACK
\cite{larsen1}, LANSO \cite{larsen,larsen1}, the MATLAB internal
function {\sf svds} and some others when computing the largest and
smallest singular triplets.

Hochstenbach \cite{Hochstenbach2001,Hochstenbach2002} shows that for
nested subspaces Ritz values approach the largest singular values
monotonically but approach the smallest ones irregularly. So the
Lanczos bidiagonalization method is more suitable for computing the
largest singular triplets and may exhibit irregular convergence
behavior when computing the smallest singular triplets. In contrast,
the smallest harmonic Ritz values converge to the smallest singular
values monotonically from above and may be better approximations. We
continue to study how to compute the smallest singular triplets more
efficiently in this paper. Based on the Lanczos bidiagonalization
process, we propose a harmonic Lanczos bidiagonalization method by
combining it with the harmonic projection principle. Our derivation
is different from that in \cite{Baglama05,Kokio04}. The method is
the same as that in \cite{Baglama05} but different from the one in
\cite{Kokio04}. We prove that for good enough projection subspaces
harmonic Ritz values converge if the columns of $A$ are strongly
linearly independent. On the other hand, harmonic Ritz values may
miss some desired singular values when the columns of $A$ are almost
linearly dependent. So harmonic Ritz values may not be reliable.
Furthermore, harmonic Ritz vectors may converge irregularly and even
may fail to converge. These results imply that either implicitly or
explicitly restarted algorithms may converge very slowly, converge
irregularly or fail to converge. To circumvent these drawbacks,
combining the harmonic projection principle with the harmonic
projection principle, we propose a refined harmonic Lanczos
bidiagonalization method that takes the Rayleigh quotients of
harmonic Ritz vectors as more accurate and reliable approximate
singular values and extracts the best approximations to the desired
singular vectors from the given subspaces that minimize the
residuals formed with the Rayleigh quotients. We prove that refined
harmonic Ritz approximations converge once the Krylov subspaces are
good enough and the Rayleigh quotients converge. We then develop an
implicitly restarted refined harmonic Lanczos bidiagonalization
algorithm (IRRHLB). Based on the refined harmonic Ritz
approximations to the desired singular vectors, in the spirit of
Jia's work \cite{jia99a,jia02}, we propose a new shifts scheme,
called the refined harmonic shifts, that we show to be theoretically
better than the harmonic shifts used within the implicitly restarted
harmonic Lanczos bidiagonalization algorithm (IRHLB) and IRLANB.
Motivated by \cite{jia99a,jia02}, we propose an efficient procedure
to compute the refined harmonic shifts accurately.

It is worth noting that Kokiopoulou {\em et al. }\cite{Kokio04}
also use the refined projection principle to compute the refined
harmonic Ritz approximations. They exploit the lower Lanczos
bidiagonalization process, use the Ritz or harmonic Ritz values as
shifts in the algorithm and compute the smallest singular triplets
one by one by exploiting deflation. They only use the refined
projection principle as refinement postprocessing at the end of
each restart. The authors demonstrate that computing refined
(harmonic) Ritz vectors and thus refined Ritz values benefits the
overall convergence process. In particular, they show that while
convergence is not apparent in terms of harmonic residual norms,
monitoring refined residuals predicts convergence more accurately
and safely. In contrast, based on the upper Lanczos
bidiagonalization process and the refined projection principle, we
propose a truly new method--the refined harmonic Lanczos
bidiagonalization method that computes refined harmonic Ritz
vectors as new approximations. We then develop IRRHLB with use of
the new better shifts, called refined harmonic shifts, based on
refined harmonic Ritz approximations. IRRHLB computes all the
desired smallest singular triplets simultaneously.

The paper is organized as follows. In \S 2, based on the Lanczos
bidiagonalization process, we derive the harmonic Lanczos
bidiagonalization method and then present some basic and important
properties of approximate singular vectors to be used later.
Exploiting Jia's results in \cite{jiaharm}, we then make a
convergence analysis. In \S 3 we propose the refined harmonic
Lanczos bidiagonalization method. We prove that the refined
harmonic Ritz approximations converge for good enough subspaces
once the Rayleigh quotients converge. In \S 4, we consider
selection of the shifts involved. For IRHLB, similar to what is
done in \cite{Baglama05,Kokio04}, we use the harmonic Ritz values.
For IRRHLB, by exploiting the available refined harmonic Ritz
approximations, we propose the refined harmonic shifts that are
proved to be theoretically better than the harmonic shifts. We
then present an efficient procedure to compute them. We show that
in finite precision the refined harmonic shifts can be computed
accurately. Meanwhile, we extend the adaptive shifting strategy
proposed by Larsen \cite{larsen1} and modified by Jia and Niu
\cite{jianiu03} to IRHLB and IRRHLB. In \S 5, we report numerical
results and compare IRRHLB with the five other  state of art
algorithms: IRHLB, IRRLB, IRLB, IRLANB and IRLBA, indicating that
IRRHLB is at least competitive with the five other  algorithms and
can be considerably more efficient when computing the smallest
singular triplets. Finally, we conclude the paper with some
remarks in \S 6.

We introduce some notations to be used. Denote by $\|\cdot\|$ the
spectral norm of a matrix and the vector 2-norm, by
$\kappa(A)=\frac{\sigma_{N}}{\sigma_1}$, by ${\cal
K}_m(C,v_1)=span\{ v_1,Cv_1,\ldots,C^{m-1}v_1 \}$ the
$m$-dimensional Krylov subspace generated by the matrix $C$ and
the starting $v_1$, by the superscript `T' the transpose of a
matrix or vector, by $I$ the identity matrix with the order clear
from the context and by $e_m$ the $m$-th coordinate vector of
dimension $m$.

\section {The harmonic Lanczos bidiagonalization method and convergence}

Golub {\em et al.} \cite{golub81} propose a Lanczos
bidiagonalization method that can compute either the largest or the
smallest singular triplets of $A$. The method is equivalent to the
symmetric Lanczos method for the eigenproblem of $\tilde A$ starting
with a special vector \cite{golub81,parlett} and is based on the
Lanczos bidiagonalization process \cite{bjorck96,golub96,paige},
which satisfies the following relations if it does not break down
before step $m$:
\begin{eqnarray}
       AQ_m&=&P_mB_m,\label{bid11} \\
       A^{\rm T}P_m&=&Q_mB_m^{\rm T}+\beta_mq_{m+1}e_m^{\rm T},\label{bid21}
    \end{eqnarray}
where the $m\times m$ matrix
\begin{equation}
    B_m=\left ( \begin{array} {c c c c}
                       \alpha_1 & \beta_1  &        &         \\
                                & \alpha_2 & \ddots &         \\
                                &          & \ddots & \beta_{m-1} \\
                                &          &        & \alpha_m
                     \end{array} \right), \label{Bm}
\end{equation}
and  the columns of $Q_m=(q_1,q_2,\ldots,q_m)$ and
$P_m=(p_1,p_2,\ldots,p_m)$ form orthonormal bases of the Krylov
subspaces ${\cal K}_m(A^{\rm T}A,q_1)$ and ${\cal K}_m(AA^{\rm
T},p_1)$, respectively.
So we have
    \begin{equation}\label{projmatrix}
       P_m^{\rm T}AQ_m=B_m.
    \end{equation}

The Lanczos bidiagonalization method computes the singular
triplets $(\tilde\sigma_i,s_i,w_i),$ $i=1,2,\ldots,m$ of $B_m$ and
then uses some of $(\tilde\sigma_i,P_ms_i,Q_mw_i)$, called the
Ritz approximations, to approximate the largest and/or smallest
singular triplets of $A$.

In finite precision, the columns of $P_m$ and of $Q_m$ may rapidly
lose orthogonality. A partial reorthogonalization strategy
\cite{larsen,larsen1} is an effective technique for maintaining
numerical orthogonality. However, Simon and Zha \cite{Simon00}
show that it generally suffices to partially reorthogonalize only
$P_m$ or $Q_m$ rather than reorthogonalizing them simultaneously.
This may reduce the computational cost considerably when only
reorthogonalizing the columns of $Q_m$ for $M\gg N$. In our codes,
we adopt the strategy from \cite{Baglama05} which is based on
\cite{Simon00}.

Given the
subspace
    $$
    E=span\left\{\left(\begin{array}{cc}
    P_m&0\\
    0&Q_m
    \end{array}
    \right)\right\},
    $$
the harmonic projection method of $\tilde{A}$ computes
$(\theta_i,\tilde\varphi_i)$ satisfying the requirements
 \begin{equation}
      \left \{ \begin{array}{rcl}
                   \tilde \varphi_i&=&\left (P_ms_i \atop Q_mw_i \right ) \in E, \\
                   (\tilde A-\theta_i I)\tilde \varphi_i&\bot&\tilde AE
               \end{array} \right. \label{proj}
    \end{equation}
and uses them as approximations to some eigenpairs of $\tilde A$
\cite{bai,stewart01,vorst}.

Making use of (\ref{projmatrix}), we see that (\ref{proj}) is
equivalent to the generalized eigenproblem
\begin{equation}
      \left ( \begin{array}{cc}
               0 & B_m \\
               B_m^{\rm T} & 0
            \end{array}\right )
      \left (s_i \atop w_i \right )=\frac{1}{\theta_i}
      \left ( \begin{array}{c c}
               B_mB_m^{\rm T}+\beta_m^2e_me_m^{\rm T} & 0 \\
               0 & B_m^{\rm T}B_m
             \end{array}\right )
      \left ( s_i \atop  w_i \right ).\label{harmbl}
    \end{equation}
$B_m$ is nonsingular as $A$ has full column rank and its singular
values interlace those of $A$ \cite[p. 449]{golub96}. This is a
symmetric positive definite generalized eigenproblem, so its
eigenvalues are all real and nonzero \cite{golub96,stewart01}.
Furthermore, we present the following result.

\begin{theorem}\label{th1}
If $(\theta_i,s_i,w_i)$ satisfies {\rm (\ref{harmbl})}, then
$(-\theta_i,s_i,-w_i)$ or equivalently $(-\theta_i,-s_i,w_i)$
satisfies {\rm (\ref{harmbl})} too, that is, the eigenvalues of
{\rm (\ref{harmbl})} come in plus-and-minus pairs and the eigenvectors
have a special structure.
\end{theorem}

{\em Proof.}  Equation (\ref{harmbl}) gives rise to
\begin{eqnarray*}
B_mw_i&=&\frac{1}{\theta_i}\left(B_mB_m^{\rm T}+\beta_m^2e_m
e_m^{\rm T}\right)s_i,\\
B_m^{\rm T}s_i&=&\frac{1}{\theta_i}B_m^{\rm T}B_mw_i.
\end{eqnarray*}
So we can readily see that the assertion holds.\qquad
\endproof

Assume that the nonnegative eigenvalues of (\ref{harmbl}) are
ordered as
$$
0\leq \theta_1\leq \theta_2\leq \cdots \leq \theta_{k+l},
$$
where $k+l=m$. Then we use
\begin{equation}\label{harmonic}
(\theta_i,\tilde u_i=P_m s_i/\|s_i\|=P_m\tilde s_i,
    \tilde v_i=Q_m w_i/\|w_i\|=Q_m\tilde w_i),i=1,2,\ldots,k
\end{equation}
as approximations to the $k$ smallest singular triplets
$(\sigma_i,u_i,v_i)$. This method is called the harmonic Lanczos
bidiagonalization method. $\theta_i$'s are called the harmonic
Ritz values, $\tilde u_i$'s and $\tilde v_i$'s the (left and
right) harmonic Ritz vectors, and $(\theta_i,\tilde u_i,\tilde
v_i)$'s the harmonic Ritz approximations. It is proved in
\cite{Hochstenbach2002} that
\begin{equation}
\sigma_i\leq\theta_i, \ i=1,2,\ldots,m.\label{mono}
\end{equation}

We have the following basic and important properties, which will
play a key role in \S 4.2.

\begin{theorem}\label{th2}
For $i\not=j$ it holds that
    \begin{equation} \label{Borth}
         s_i^{\rm T}B_mw_j=0,\, w_i^{\rm T}B_m^{\rm T}s_j=0
    \end{equation}
and
\begin{equation}\label{Borth2}
         \tilde u_i^{\rm T}A\tilde v_j=0,\,
         \tilde v_i^{\rm T}A^{\rm T}\tilde u_j=0.
    \end{equation}
\end{theorem}

{\em Proof.} Since (\ref{harmbl}) is a symmetric positive definite
generalized eigenproblem, for $i\not= j$ we have
    $$
        (s_i^{\rm T},w_i^{\rm T}) \left ( \begin{array}{cc}
               0 & B_m \\
               B_m^{\rm T} & 0
            \end{array}\right )\left (s_j \atop w_j \right ) =0,
    $$
from which it follows that (\ref{Borth}) holds.

(\ref{Borth2}) follows from $\tilde u_i=P_ms_i/\|s_i\|,\,\tilde
v_i=Q_mw_i/\|w_i\|$ and (\ref{projmatrix}). \qquad
\endproof

We now discuss efficient computation of harmonic Ritz
approximations. It is disappointing that (\ref{harmbl}) is a
$2m\times 2m$ generalized eigenproblem. Fortunately, we can reduce
(\ref{harmbl}) to a half sized SVD problem that can be solved more
cheaply and accurately, as shown below.

We get from (\ref{harmbl})
\begin{eqnarray}
(B_mB_m^{\rm T}+\beta_m^2e_me_m^{\rm T})s_i&=&\theta_i B_mw_i,\label{expand1}\\
B_mw_i&=&\theta_is_i,\label{expand2}
\end{eqnarray}
from which it follows that
\begin{equation}
 (B_mB_m^{\rm T}+\beta_m^2e_me_m^{\rm
T})s_i=\theta_i^2s_i. \label{reduction}
\end{equation}
So, the $(\theta_i,s_i)$'s are the singular values and right
singular vectors of the $(m+1)\times m$ matrix
\begin{equation} \label{svd1}
\left(\begin{array}{c}
B_m^{\rm T}\\
\beta_me_m^{\rm T}
\end{array}
\right)
\end{equation}
and the left singular vectors
\begin{equation}
w_i=\theta_iB_m^{-1}s_i. \label{bidiag}
\end{equation}
Therefore, we can get $(\theta_i,s_i)$ more accurately and
efficiently by computing the SVD of the half sized (\ref{svd1}) and
obtain all the $w_i$'s by solving the bidiagonal linear systems
$B_mw_i=\theta_is_i$ at a total cost of $O(m^2)$ flops.

We comment that, based on the harmonic projection of $A^{\rm T}A$
onto ${\cal K}_m(A^{\rm T}A,q_1)$, Baglama and Reichel
\cite{Baglama05,baglama06} also derive
(\ref{reduction})--(\ref{bidiag}). Our method is the same as that in
\cite{Baglama05} but is different from the one in \cite{Kokio04},
which is based on the lower Lanczos bidiagonalization process.

From (\ref{bid11}) and (\ref{bid21}), we have
    $$\|A\tilde v_i-\theta_i\tilde u_i\|=\|AQ_m\tilde w_i-\theta_iP_m\tilde s_i\|
      =\|P_mB_m\tilde w_i-\theta_iP_m\tilde s_i\|=\|B_m\tilde w_i-\theta_i\tilde s_i\|$$
and
$$
    \|A^{\rm T}\tilde u_i-\theta_i\tilde v_i\|=\sqrt{\|B_m^{\rm T}
    \tilde s_i-\theta_i\tilde w_i\|^2+
      \beta_m^2|e_m^{\rm T}\tilde s_i|^2}.
$$
Therefore, if
\begin{eqnarray}
\sqrt{\|A\tilde v_i-\theta_i\tilde u_i\|^2+\|A^{\rm T}\tilde
        u_i-\theta_i\tilde v_i\|^2}&=&\sqrt{\|B_m\tilde w_i-\theta_i\tilde s_i\|^2+
        \|B_m^{\rm T}\tilde s_i-\theta_i \tilde w_i\|^2+
      \beta_m^2|e_m^{\rm T}\tilde s_i|^2}\nonumber\\
      &<& tol, \label{harmstop}
\end{eqnarray}
where $tol$ is a user prescribed accuracy, then the method is
accepted as converged for $tol$.

We now analyze the convergence. Jia \cite{jiaharm} establishes a
general convergence theory of harmonic projection methods for large
eigenproblems. The theory can be adapted here.

From (\ref{harmbl}), set the matrices
$$
\tilde B=\left ( \begin{array} {c c}
                      0     & B_m  \\
                      B_m^{\rm T} & 0
                   \end{array} \right )
$$
and
$$\tilde C=\left ( \begin{array}{c c}
               B_mB_m^{\rm T}+\beta_m^2e_me_m^{\rm T} & 0 \\
               0 & B_m^{\rm T}B_m
             \end{array}\right ).
$$
Recall from Theorem~\ref{th1} that $\pm\theta_i$'s are the
eigenvalues of $\tilde B^{-1}\tilde C$. The following result is
direct from Theorem 2.1 and Corollary 2.2 of \cite{jiaharm}.

\begin{theorem}\label{theo1}
Assume that $(\sigma,u,v)$ is a singular triplet of $A$ and define
$\varepsilon=\sin \angle \left ( \left (
    u \atop v \right ),E \right )$ to be the distance between the vector
$\left(
    u \atop v \right )$ and the subspace $E$. Then there exists a perturbation
    matrix $F$ satisfying
    \begin{equation}
       \|F\|\leq \frac{\varepsilon}{\sqrt{1-\varepsilon^2}}\|\tilde B^{-1}\|(
       \sigma\|A\|+\|A\|^2),\label{value2}
    \end{equation}
such that the exact singular value $\sigma$ of $A$ is an
eigenvalue of $\tilde B^{-1}\tilde C+F$. Furthermore, there exists
an eigenvalue $\theta$ of $\tilde B^{-1}
    \tilde C$ satisfying
    \begin{equation}
       |\theta-\sigma|\leq (2\|A\|+\|F\|)\|F\|.
    \end{equation}
\end{theorem}

This theorem shows that if $\varepsilon$ tends to zero and $\|\tilde
B^{-1}\|$ is uniformly bounded  then there always exists one
harmonic Ritz value $\theta$ that converges to the desired singular
value $\sigma$. The interlacing theorem of singular values \cite[p.
449]{golub96} tells us that
    $$
    \frac{1}{\sigma_N}\leq \|\tilde B^{-1}\|=\|B_m^{-1}\|\leq
    \frac{1}{\sigma_1}=\|A^+\|
    $$
is uniformly bounded. As a result, if $\varepsilon=\sin \angle \left
( \left (
    u \atop v \right ),E \right )\rightarrow 0$ʱ, we should have
$\theta\rightarrow\sigma$.

However, the situation is by no means so simple, and it is
instructive to see what will happen when only speaking of
$\varepsilon$ small. Note that
$$
\|B_m^{-1}\|\|A\|\leq\|A^+\|\|A\|=\kappa(A).
$$
So if $\kappa(A)=O(\frac{1}{\varepsilon})$, that is, the columns
of $A$ are almost linearly dependent, then $\|F\|$ may not be near
zero, so that $|\theta-\sigma|$ may not be small. Actually,
$\|B_m^{-1}\|\rightarrow \frac{1}{\sigma_1}=\|A^+\|$ once the
smallest singular value (Ritz value) of $B_m$ converges to
$\sigma_1$. In this case, the harmonic Lanczos bidiagonalization
method may miss $\sigma$. So the method may not be reliable, and
$\theta$ is only guaranteed to be a good approximation to $\sigma$
only if $\varepsilon$ is very small and the columns of $A$ are
strongly linearly independent.

To improve convergence and reliability of the method, we recommend
the Rayleigh quotient $\rho_i=\tilde u_i^TA\tilde v_i=\tilde
s_i^TB_m\tilde w_i$ as a new approximation to $\sigma_i$, as was
also done in \cite{Hochstenbach2002}. $\rho_i$ is more accurate and
reliable than $\theta_i$. Correspondingly, $\theta_i$ in
(\ref{harmstop}) is replaced by $\rho_i$. We refer to \cite{jiaharm}
for more theoretical results and arguments on such a replacement.

The following result is a direct application of Theorem 3.2 of
\cite{jiaharm}.

\begin{theorem}\label{theo2}
Let $(\theta,z=\left ( s \atop w \right ))$ be an eigenpair of
$\tilde B^{-1}\tilde C$, and assume that $Z_\bot$ is such that the
square matrix $\left ( z,Z_\bot \right )$ is orthogonal and
transforms $\tilde B^{-1}\tilde C$ into
    \begin{equation}
       \left ( z^{\rm T} \atop Z_\bot^{\rm T} \right )
       \tilde B^{-1}\tilde C \left (z,Z_\bot \right )
       =\left ( \begin{array} {c c}
                   \theta & g^{\rm T} \\
                   0      & G
                \end{array} \right ).
    \label{decomp2}
    \end{equation}
Then if
    \begin{equation}
       {\rm  sep}(\sigma,G)>0,
       \label{uniform2}
    \end{equation}
it holds that
    \begin{eqnarray}
        \sin \angle
        \left (\left ( u\atop v \right),\left (\tilde u \atop \tilde v\right )\right )
        &\leq& \left ( 1+ \frac {2\|\tilde B^{-1}\|\|A\|} {\sqrt {1-\varepsilon^2}
        {\rm  sep}(\sigma,G)} \right )\varepsilon\label{vector2}\\
        &\leq& \left ( 1+ \frac {2\|\tilde B^{-1}\|\|A\|} {\sqrt {1-\varepsilon^2}
        ({\rm  sep}(\theta,G)-|\sigma-\theta|)}\right)\varepsilon. \label{sep}
    \end{eqnarray}
\end{theorem}

Theorems~\ref{theo2} states that if the separation ${\rm
sep}(\theta,G)$ of $\theta$ and the other harmonic Ritz values is
uniformly bounded below by a positive constant then the harmonic
Ritz approximations $\tilde{u},\tilde{v}$ converge. Unfortunately,
however, for a general $A$, ${\rm sep}(\theta,G)$ can be
arbitrarily near zero, i.e., $\theta$ can be arbitrarily close to
the other harmonic Ritz values. As a result, the upper bounds
(\ref{vector2}) and (\ref{sep}) can converge to zero very slowly
and irregularly and even fail to do so as $\varepsilon\rightarrow
0$. This means that the approximate singular vectors may converge
very slowly and irregularly and even may fail to converge.

\section{The refined harmonic Lanczos bidiagonalization method}

The previous analysis shows that the harmonic Ritz approximations
may converge slowly and irregularly and even may fail to converge.
To overcome this intrinsic drawback, we now combine the harmonic
Lanczos bidiagonalization method with the refined projection
principle and derive a refined harmonic Lanczos bidiagonalization
method. Recall that $\rho_i=\tilde s_i^TB_m\tilde w_i$. We use
$(\rho_i,\tilde \psi_i)$ satisfying
\begin{equation}
      \left \{ \begin{array}{rcl}
                   \tilde \varphi_i&=&\left (P_ms_i \atop Q_mw_i \right ) \in E, \\
                   (\tilde A-\theta_i I)\tilde \varphi_i &\bot& \tilde AE,\\
                   \|\tilde A \tilde \psi_i-\rho_i \tilde
                   \psi_i\|&=&
                      \min_{\psi \in E,\|\psi\|=1}\|\tilde A\psi-\rho_i \psi\|
               \end{array} \right .
               \label{def}
\end{equation}
 to replace
 $(\theta_i,\tilde \varphi_i)$ as a new approximation to an eigenpair of
 $\tilde A$. The method first uses the harmonic Lanczos bidiagonalization method to
 compute $(\theta_i,\tilde\varphi_i)$
 and then forms the Rayleigh quotient $\rho_i=\tilde
s_i^TB_m\tilde w_i$. With each $\rho_i,\,i=1,2,\ldots,k$, it
computes $\psi_i$ by solving the minimization problem.

The following results adapted from \cite{jia97,jia02} can be used
to compute $\tilde{\psi}_i$ efficiently and accurately.

\begin{theorem}\label{theo4}
Let $z_i=(x_i^{\rm T},y_i^{\rm T})^{\rm T}$ be the right singular
vector of the matrix
  $$
  \left (
     \begin{array} {c c}
        0                & B_m \\
        B_m^{\rm T}            & 0 \\
        \beta_me_m^{\rm T} & 0
      \end{array} \right )
  -\rho_i \left (
      \begin{array} {c c}
           I & 0 \\
           0 & I \\
           0 & 0
      \end{array} \right )
  $$
associated with its smallest singular value $\sigma_{\min}$. Then
  \begin{equation}
     \tilde \psi_i=\left ( \begin{array} {c c}
                             P_m & 0  \\
                             0   & Q_m
                             \end{array} \right ) z_i,
  \end{equation}
  \begin{equation}
 \|\tilde A \tilde \psi_i- \rho_i \tilde \psi_i\|=\sigma_{\min}.
  \end{equation}
\end{theorem}

With $\tilde{\psi}_i$ at hand, we define the new left and right
approximate singular vectors as
  \begin{equation}
  \hat u_i=P_mx_i/\|x_i\|=P_m\hat x_i,\, \hat v_i=Q_my_i/\|y_i\|=Q_m\hat
  y_i
  \label{refinevector}
  \end{equation}
and use $(\rho_i,\hat u_i,\hat v_i)$'s to approximate the $k$
smallest singular triplets of $A$. We call $(\rho_i,\hat u_i,\hat
v_i)$ a refined harmonic Ritz triplet and $\hat u_i,\hat v_i$ a
refined harmonic Ritz approximation.

Similar to (\ref{harmstop}), $(\rho_i,\hat u_i,\hat v_i)$ is
accepted as converged if
\begin{eqnarray}
\sqrt{\|A\hat v_i-\rho_i\hat u_i\|^2+\|A^{\rm T}\hat v_i-\rho_i\hat
  v_i\|^2}&=&\sqrt{\|B_m\hat y_i-\rho_i\hat x_i\|^2+\|B_m^{\rm T}\hat
x_i-\rho_i\hat y_i\|^2
   +\beta_m^2|e_m\hat x_i|^2}\nonumber\\
   &\leq& tol.\label{refinedstop}
\end{eqnarray}

The following result is taken directly from Theorem 4.1 of
\cite{jiastew}.

\begin{theorem}\label{theo5}
Let $(\sigma,u,v)$ be a singular triplet of $A$. Then there exist
$U_\bot$ and $V_\bot$ such that $(u,U_\bot)$ and $(v,V_\bot)$
are orthogonal and
\begin{equation}
     \left ( u^{\rm T} \atop U_\bot^{\rm T} \right )A(v,V_\bot)=
     \left ( \begin{array} {c c}
                \sigma & 0 \\
                0      & L
             \end{array}  \right ),
  \end{equation}
where $L=U_\bot^{\rm T}AV_\bot$. Set
$$
\tilde L=
  \left ( \begin{array} {c c}
             0    & L \\
             L^{\rm T}  & 0
          \end{array}  \right )
$$
and assume that $(\rho,\hat u,\hat v)$ is approximating
$(\sigma,u,v)$. Then if
  \begin{equation}
     {\rm  sep}(\rho,\tilde L)\geq {\rm  sep}(\sigma,\tilde L)-|
     \rho- \sigma|>0,
  \end{equation}
we have
  \begin{equation}
     \sin\theta\left ( \left (\hat u \atop \hat v \right ),
     \left (u \atop v\right ) \right )
     \leq \frac {\|\tilde A-\rho I\|\varepsilon+
     |\rho-\sigma|}
                {\sqrt {1-\varepsilon^2}({\rm  sep}(\sigma,\tilde L)-
                |\rho-\sigma|)}.
  \end{equation}
\end{theorem}

Note that ${\rm sep}(\sigma,\tilde L)$ is the gap of $\sigma$ and
the other singular values of $A$ and is a fixed constant.
Theorems~\ref{theo5} shows that the refined harmonic Ritz
approximations converge once $\varepsilon\rightarrow 0$ and
$\rho\rightarrow\sigma$. Therefore, the refined harmonic Lanczos
bidiagonalization method overcomes, to great extent, the possible
non-convergence of the harmonic Ritz approximations. 

\section{The implicit restarting technique, shifts selection and an adaptive
shifting strategy}

\subsection{The implicit restarting technique}

Due to the storage requirements and computational cost, in
practice, the number of steps cannot be large and must be limited.
For a relatively small $m$, however, the $m$-dimensional subspaces
${\cal K}_m(A^{\rm T}A,q_1)$ and ${\cal K}_m(AA^{\rm T},p_1)$, in
general, do not contain enough information on the desired right
and left singular vectors, so that both the harmonic and refined
harmonic Lanczos bidiagonalization methods do not converge.
Therefore, it is necessary to restart the methods. The idea is to
repeatedly update new starting vectors based on the information
available and construct increasingly better Krylov subspaces until
the methods converge. Implicit restarting is usually preferable
not only because of efficiency of the restart procedure, but also
because the implicit procedure is more effective at locking in
desired directions and purging unwanted ones.

We briefly review the implicit restarting technique for the Lanczos
bidiagonalization process \cite{bjorck94,larsen1}. Note that
$m=k+l$. After running $l$ implicit QR iteration steps on $B_m$
using the shifts $\mu_j,j=1,2,\ldots,l$, we get
\begin{equation}
    \left \{ \begin{array} {l}
               (B_m^{\rm T}B_m-\mu_1^2 I)\cdots (B_m^{\rm T}B_m-\mu_l^2 I)
               =\tilde P R, \\
               {\tilde P}^{\rm T}B_m\tilde Q\mbox {\ \ \ \ upper
               bidiagonal},
             \end{array} \right.
  \end{equation}
where $\tilde P$ and $\tilde Q$ are the accumulations of Givens
rotations applied to $B_m$ from the left and right, respectively.
Define $Q_m^+=Q_m\tilde Q,\,P_m^+=P_m\tilde P$ and $B_m^+={\tilde
P}^{\rm T}B_m\tilde Q$. This process is achieved implicitly from
$B_m^{\rm T}B_m$ to $(B_m^+)^{\rm T}B_m^+$ by working on $B_m$
directly.

Performing the above $l$ implicit QR iteration steps gives the following
relations \cite{bjorck94}:
  \begin{eqnarray}
  AQ_k^+&=&P_k^+B_k^+,\\
  A^{\rm T}P_k^+&=&Q_k^+{B_{k}^+}^{\rm  T}+(\beta_k\tilde p_{m,k}
  q_{m+1}+\beta_k^+ q_{k+1}^+)e_k^{\rm  T},
  \end{eqnarray}
where $\tilde p_{m,k}$ is the entry of $\tilde P$ in position
$(m,k)$ and the updated starting vector has the form
\begin{equation}
\gamma q_1^+=\prod_{j=1}^l(A^{\rm T}A-\mu_j^2 I)q_1 \label{update}
\end{equation}
with $\gamma$ a factor making $\|q_1^+\|=1$. Since $\beta_k\tilde
p_{m,k}q_{m+1}+\beta_k^+ q_{k+1}^+$ is orthogonal to $Q_{k}^+$, we
have obtained a $k$-step Lanczos bidiagonalization process starting
with $q_1^+$. It is then extended to a $m$-step Lanczos
bidiagonalization process in a standard way. So we avoid restarting
the process from scratch and do it from step $k+1$ upwards. This
saves the computational cost of the first $k$ steps of the process.
Applying the implicit restarting technique to the harmonic Lanczos
bidiagonalization method and its refined version in such a way, we
have formally sketched an implicitly restarted harmonic Lanczos
bidiagonalization algorithm (IRHLB) and an implicitly restarted
refined harmonic Lanczos bidiagonalization algorithm (IRRHLB).

\subsection{Shifts selection}

We can run IRHLB and IRRHLB once the shifts $\mu_1,\mu_2,\\
\ldots,\mu_l$ are given. However, in order to make them work as
efficiently as possible, we should select the best possible shifts
in some sense for each algorithm. In the same spirit of
\cite{jia99a,jia02}, it has been shown in \cite{jianiu03} that if
the shifts are more accurate approximations to some of the
unwanted singular values of $A$ then the resulting subspaces
contain more information on the desired singular vectors. The
better the subspaces are, the faster IRHLB and IRRHLB may
converge. For eigenproblems and SVD problems, Morgan
\cite{Morgan00,Morgan06} and Kokiopoulou {\em et al}.
\cite{Kokio04} suggest using unwanted harmonic Ritz values as
shifts, called the harmonic shifts here. These shifts are natural
choices as they are the best approximations available to some of
the unwanted eigenvalues and the unwanted singular values,
respectively. So, for our IRHLB we also use the $l$ unwanted
harmonic Ritz values $\theta_{k+1},\theta_{k+2},\ldots,\theta_m$
as shifts. Since the refined harmonic approximations $\hat
u_i,\hat v_i$ are optimal in the sense of residual minimizations,
they are generally more accurate than the harmonic Ritz
approximations $\tilde u_i,\tilde v_i$. Therefore, based on $\hat
u_i,\hat v_i,\,i=1,2,\ldots,k$, it should be possible to find
better possible shifts than the harmonic shifts.

The following important result on the harmonic shifts is crucial for
us to introduce and understand new better shifts for use within
IRRHLB.

\begin{theorem}\label{theo6}
Define
\begin{eqnarray*}
&\tilde U=(\tilde u_1,\ldots,\tilde u_k),\ \tilde V=(\tilde
v_1,\ldots,\tilde v_k),& \\
&\tilde{U}_\perp=(\pm\tilde u_{k+1},\ldots,\pm\tilde u_m),\
  \tilde{V}_\perp=(\pm\tilde v_{k+1},\ldots,\pm\tilde v_m).&
\end{eqnarray*}
Then the harmonic shifts
$\theta_{k+1},\theta_{k+2},\ldots,\theta_m$ are the absolute
values of the $l$ harmonic Ritz values of $\tilde A$ with respect
to the subspace $span\{(\tilde{U}_\perp^{\rm
T},\tilde{V}_\perp^{\rm T})^{\rm T}\}$.
\end{theorem}

{\em Proof.} From definition (\ref{proj}) of the harmonic projection
as well as the relationship between (\ref{proj}) and (\ref{harmbl}),
it is easily verified that for $i=k+1,\ldots,m$, if the $i$-th
column of $\tilde{U}_\perp$ and that of $\tilde{V}_\perp$ have the
same or opposite $\pm$ sign, then $\tilde A$ has $\theta_i$ or
$-\theta_i$ as one harmonic Ritz value with respect to the subspace
$span\{(\tilde{U}_\perp^{\rm T},\tilde{V}_\perp^{\rm T})^{\rm
T}\}$.\qquad
\endproof

We see from (\ref{Borth2}) that
$$
\tilde U^{\rm
T}A\tilde{V}_\perp=0,\,\tilde V^{\rm T}A^{\rm T}\tilde{U}_\perp=0
$$
and
$$
  span\{\tilde V\}\oplus span\{\tilde{V}_\perp\}=span\{Q_m\},\,
  span\{\tilde U\}\oplus span\{\tilde {U}_\perp\}=span\{P_m\}.
$$

Define $\hat U=(\hat u_1,\ldots,\hat u_k)$ and $\hat V=(\hat
v_1,\ldots,\hat v_k)$, and let $\hat{U}_\perp,\hat{V}_\perp$ be
matrices with $l=m-k$ columns satisfying
\begin{equation}
 \hat{U}^{\rm
T}A\hat{V}_\perp=0,\hat{V}^{\rm T}A^{\rm T}
  \hat{U}_\perp=0 \label{Aorth}
\end{equation}
and
  \begin{equation}
  span\{\hat V\}\oplus span\{\hat{V}_\perp\}=span\{Q_m\},\,
  span\{\hat U\}\oplus
  span\{\hat{U}_\perp\}=span\{P_m\},\label{sumdecomp}
  \end{equation}
where $\oplus$ denotes the direct sum. Jia \cite{jia04} derives a
number of theoretical results that compare refined Ritz vectors
and Ritz vectors. At this moment, we temporarily regard $\tilde A$
as a general matrix, and $\tilde\varphi_i,\tilde\psi_i$ are a Ritz
and the corresponding refined Ritz vector of $\tilde A$ with
respect to a general subspace $E$, respectively. One of Jia's
results says that we always have
$$
\|(\tilde A-\rho_i I)\tilde\psi_i\|<\|(\tilde A-\rho_i
I)\tilde\varphi_i\|
$$
if the left-hand side is not zero and
$$
\|(\tilde A-\rho_i I)\tilde\psi_i\|\ll\|(\tilde A-\rho_i
I)\tilde\varphi_i\|
$$
may occur if $\rho_i$ is close to some $\theta_j$ for $j\not=i$.
By standard perturbation theory in terms of residual norms, these
two results demonstrate that $\tilde\psi_i$ is more accurate and
can be much more accurate than $\tilde\varphi_i$. Here we should
point out that these claims hold without requiring that $E$ is
sufficiently good. Jia \cite{jia04} constructs a number of
symmetric matrices having well separated simple eigenvalues and
accurate subspaces to illustrate this. More precisely, assuming
that $\rho_i$ and $\tilde\varphi_i,\,\tilde\psi_i$ are used to
approximate the eigenvalue $\lambda_i$ and the eigenvector
$\varphi_i$ of $\tilde A$ and $\varepsilon$ is the distance
between $\varphi_i$ and the subspace $E$, Jia's examples show that
we can indeed have
$$
|\rho_i-\lambda_i|=O(\varepsilon),\ \|(\tilde A-\rho_i
I)\tilde\psi_i\|=O(\varepsilon), \ \|(\tilde A-\rho_i
I)\tilde\varphi_i\|=O(1).
$$
The above results are easily adapted to the harmonic and refined
harmonic Ritz vectors. Some similar symmetric matrices are
constructed by Jia in \cite{jiaharm} for which the harmonic Ritz
vectors have no accuracy at all but refined harmonic ones have
accuracy $O(\varepsilon)$. Coming back to our SVD context, the
above results and analysis indicate that $\hat u_i$ and $\hat v_i$
are more accurate and can be much more accurate than $\tilde u_i$
and $\tilde v_i$ without the assumption that projection subspaces
are sufficiently good.

Based on the above, it is evident that the subspaces
$span\{\hat{U}_\perp\}$ and $span\{\hat{V}_\perp\}$ contain
(possibly much) more accurate approximations to $\pm u_i$ and $\pm
v_i$, $i=k+1,k+2,\ldots,N$ than the subspaces
$span\{\tilde{U}_\perp\}$ and $span\{\tilde{V}_\perp\}$ do. This, in
turn, means that the subspace $span\{(\hat{U}_\perp^{\rm
T},\hat{V}_\perp^{\rm T})^{\rm T}\}$ contains (possibly much) more
accurate approximations to the eigenvectors
$\frac{1}{\sqrt{2}}(u_i^{\rm T},\pm v_i^{\rm T})^{\rm T}$ associated
with the eigenvalues $\pm \sigma_i$, $i=k+1,k+2,\ldots,N$, of
$\tilde A$ than the subspace $span\{(\tilde{U}_\perp^{\rm
T},\tilde{V}_\perp^{\rm T})^{\rm T}\}$ does. Recall from
Theorem~\ref{theo1} that a better subspace should generally produce
more accurate harmonic Ritz values. Hence, combining with
Theorem~\ref{theo6}, we have come to the following key result.

\begin{theorem}\label{theo9}
As approximations to some of $\sigma_{k+1},\ldots,\sigma_N$, the
absolute values of the harmonic Ritz values $\xi_i,\,i=1,2,\ldots,l$
of $\tilde A$ with respect to the subspace
$span\{(\hat{U}_\perp^{\rm T},\hat{V}_\perp^{\rm T})^{\rm T}\}$ are
more accurate and can be much more accurate than the harmonic shifts
$\theta_{k+1},\theta_{k+2},\ldots,\theta_m$.
\end{theorem}

This theorem holds without assuming that $span\{Q_m\}$ and
$span\{P_m\}$ are sufficiently good. It suggests that we use
better $|\xi_i|,\,i=1,2,\ldots,l$ as shifts for use within IRRHLB.
We call them the refined harmonic shifts.

Computationally, at first glance, it seems quite complicated and
expensive to get the refined harmonic shifts as it involves
constructing $\hat{U}_\perp,\hat{V}_\perp$ that are related with
the large $\tilde A$. Inspired by the tricks in
\cite{jia99a,jia02}, however, we can exploit (\ref{projmatrix}) to
propose an efficient procedure for computing them accurately, as
shown below.

Recall (\ref{refinevector}) and define $\hat X=(\hat
x_1,\ldots,\hat x_k)$ and $\hat Y=(\hat y_1,\ldots,\hat y_k)$. We
use Householder transformations to compute the full QR
decompositions
\begin{equation}
B_m^{\rm T}\hat X=Q_{X}\left(\begin{array}{c}
R_X\\
0
\end{array}
\right),\ \  B_m\hat Y=Q_{Y}\left(\begin{array}{c} R_Y\\
0\end{array} \right), \label{qr}
\end{equation}
which costs $O(m^3)$ flops. Partition
\begin{equation}
Q_X=(Q_{X1},Q_{X2}),\ Q_Y=(Q_{Y1},Q_{Y2}),\label{partition}
\end{equation}
where $Q_{X1}$ and $Q_{Y1}$ are the first $k$ columns of $Q_X$ and
$Q_Y$, respectively, and let
\begin{equation}
\hat{U}_\perp=P_mQ_{Y2},\ \hat{V}_\perp=Q_mQ_{X2}. \label{perp}
\end{equation}
Then it can be readily verified that
\begin{eqnarray*}
 & \hat{U}^{\rm T}A\hat{V}_\perp=\hat X^{\rm T}P_m^{\rm T}AQ_mQ_{X2}
  =\hat X^{\rm T}B_mQ_{X2}=0,&\\
& \hat{V}^{\rm T}A^{\rm T}\hat{U}_\perp=\hat Y^{\rm T}Q_m^{\rm
T}A^{\rm T}P_mQ_{Y2}
  =\hat Y^{\rm T}B_m^{\rm T}Q_{Y2}=0.&
\end{eqnarray*}
So $\hat U$ and $\hat V$ defined in this way meet conditions
(\ref{Aorth}) and (\ref{sumdecomp}) and are just what we need. By
(\ref{sumdecomp}), we have
\begin{equation} span\{\hat
U\}=span\{P_mQ_{Y1}\},\ span\{\hat
V\}=span\{Q_mQ_{X1}\}.\label{newbasis}
\end{equation}

The harmonic Ritz values $\xi_i,\,i=1,2,\ldots,l$ of $\tilde A$
with respect to $(\hat{U}_\perp^{\rm T},\hat{V}_\perp^{\rm
T})^{\rm T}$ satisfy
  $$
  (\hat{U}_\perp^{\rm T},\hat{V}_\perp^{\rm T})\tilde A
  \left (\hat{U}_\perp \atop \hat{V}_\perp
  \right )g_i
    =\frac{1}{\xi_i}(\hat{U}_\perp^{\rm T},\hat{V}_\perp^{\rm T})
    \tilde A^{\rm T}\tilde A\left (
    \hat{U}_\perp \atop \hat{V}_\perp
    \right )g_i.
  $$
Exploiting (\ref{projmatrix}), we get a $l\times l$ symmetric
positive definite generalized eigenvalue problem
 \begin{eqnarray} \label{shift}
      &(Q_{Y2}^{\rm T},Q_{X2}^{\rm T})
      \left ( \begin{array}{cc}
               0 & B_m \\
               B_m^{\rm T} & 0
            \end{array}\right )
      \left (Q_{Y2} \atop Q_{X2} \right )g_i=&\nonumber\\
      &\frac{1}{\xi_i}
      (Q_{Y2}^{\rm T},Q_{X2}^{\rm T})
      \left ( \begin{array}{c c}
               B_mB_m^{\rm T}+\beta_m^2e_me_m^{\rm T} & 0 \\
               0 & B_m^{\rm T}B_m
             \end{array}\right )
      \left (Q_{Y2} \atop Q_{X2} \right )g_i.&
    \end{eqnarray}
The $\xi_i$'s are computed by the QZ algorithm
\cite{golub96,stewart} using $O(l^3)$ flops. So the total cost of
computing the refined harmonic shifts is $O(m^3)$ flops, negligible
compared with the harmonic Lanczos bidiagonalization method.

We give more details on computation of the refined harmonic
shifts. Denote by $F_m$ and $G_m$ the matrices of the left and
right-hand sides in (\ref{shift}), respectively, and observe that
{\small
\begin{eqnarray}
&F_m=Q_{Y2}^{\rm T}B_mQ_{X2}+Q_{X2}^{\rm T}B_m^{\rm
T}Q_{Y2},&\label{fm}\\
&G_m=(B_m^{\rm T}Q_{Y2})^{\rm T}(B_m^{\rm
T}Q_{Y2})+\beta_m^2(e_m^{\rm T}Q_{Y2})^{\rm T}(e_m^{\rm
T}Q_{Y2})+(B_mQ_{X2})^{\rm T}(B_mQ_{X2}).&\label{gm}
\end{eqnarray}}
Noting that the two matrices in $F_m$ are transposes each other,
we only need to form $Q_{Y2}^{\rm T}B_mQ_{X2}$ by computing
$(Q_{Y2}^{\rm T}B_m)Q_{X2}$ or $Q_{Y2}^{\rm T}(B_mQ_{X2})$, and
$Q_{Y2}^{\rm T}B_m$ or $B_mQ_{X2}$ is then used to form $G_m$.
Since $G_m$ is symmetric, we only need to compute its upper
triangular part. The total cost of forming $F_m$ and $G_m$ is
$O(m^3)$ flops. We then compute the eigenvalues
$\frac{1}{\xi_i}$'s of the symmetric positive definite matrix
pencil $(F_m,G_m)$.

Now we show that in finite precision the above procedure is
numerically stable and can compute the refined harmonic shifts
accurately. There are three major steps in the procedure: the QR
decompositions in (\ref{qr}), computation of $F_m$ and $G_m$ and
the solution of the eigenvalue problem of $(F_m, G_m)$ by the QZ
algorithm. Note that the QR decompositions can be computed using
Householder transformations in a numerically stable way (we use
the Matlab built-in code {\sf qr} in our implementation) and the
QZ algorithm are numerically stable. Therefore, omitting details
on roundoff errors, we finally compute the eigenvalues
$\frac{1}{\tilde \xi_i}$'s of a perturbed matrix pencil
$(F_m+\delta F_m,G_m+\delta G_m)$, where $\delta F_m$ and $\delta
G_m$ are the matrices of roundoff error accumulations and satisfy
\begin{equation}\label{backerror}
\frac{\|\delta F_m\|_F}{\|F_m\|_F},\,\frac{\|\delta
G_m\|_F}{\|G_m\|_F}=O(\epsilon_{\rm mach})
\end{equation}
with $\epsilon_{\rm mach}$ being the machine precision and
$\|\cdot\|_F$ the Frobenius norm.

For an eigenpair $(\frac{1}{\xi},g)$ of the pencil $(F_m,G_m)$
with $\|g\|=1$, let $\alpha=g^{\rm T}F_mg$ and $\beta=g^{\rm
T}G_mg$, so that $(\beta,\alpha)$ is a projective representation
of the eigenvalue $\frac{1}{\xi}$ \cite[p. 135]{stewart01}. Then
it is known \cite[p. 233]{stewart01} that there is an eigenvalue
$\frac{1}{\tilde{\xi}}$ of the matrix pencil $(F_m+\delta
F_m,G_m+\delta G_m)$ such that the chordal distance
\begin{equation}
\chi(\frac{1}{\xi},\frac{1}{\tilde\xi})\leq\frac{\sqrt{\|\delta
F_m\|_F^2+\|\delta
G_m\|_F^2}}{\sqrt{\alpha^2+\beta^2}}+O(\epsilon_{\rm
mach}^2).\label{chi}
\end{equation}

It is important to point out that for not too small $\xi$ the choral
distance behaves like the ordinary distance $\mid\xi-\tilde\xi\mid$;
see a remark in \cite[p. 140]{stewart01}. So, how accurate
$\tilde\xi$ is depends on the $\frac{1}{\xi}$'s condition number
\begin{equation}
\nu=\frac{1}{\sqrt{\alpha^2+\beta^2}}.\label{condnum}
\end{equation}
If one of $\alpha$ and $\beta$ is not small, $\nu$ is not large and
thus by (\ref{backerror}) the relative error of $\tilde\xi$ is
$O(\epsilon_{\rm mach})$ if $|\xi|$ is not very small.

We look at the smallest $|\xi|$. By Theorem~\ref{theo9}, it is
known that the absolute values $|\xi|$'s better approximate some
of $\sigma_{k+1},\ldots,\sigma_N$ than
$\theta_{k+1},\ldots,\theta_m$. Furthermore, recall from
(\ref{mono}) that $\theta_{k+1}\geq\sigma_{k+1}$. So the smallest
$|\xi|$ is approximately bounded below by $\sigma_{k+1}$.

In the following, we establish lower bounds for $|\alpha|$ and
$|\beta|$ and an upper bound for $\nu$ rigorously and prove when
our proposed procedure can numerically compute the refined
harmonic Ritz shifts accurately.

\begin{theorem}\label{theo10}
For any refined harmonic Ritz shift $\xi$, we have
\begin{equation}
|\alpha|\geq 2\sigma_{k+1} \mbox{ and }|\beta|\geq 2\sigma_{k+1}^2
\label{alphabound}
\end{equation}
and the $\frac{1}{\xi}$'s condition number is bounded from above:
\begin{equation}
\nu\leq\frac{1}{2\sigma_{k+1}\sqrt{1+\sigma_{k+1}^2}}.
\label{nubound}
\end{equation}
If $\sigma_{k+1}$ is not very small, then numerically the procedure
described can compute the refined harmonic shifts $|\xi_i|$'s with
relative accuracy $O(\epsilon_{\rm mach})$.
\end{theorem}

{\em Proof}. We prove (\ref{nubound}) and (\ref{xibound}) in turn.
To prove (\ref{nubound}), we first estimate $|\alpha|$ and then
$|\beta|$. From the definition of $F_m$ and $G_m$, we have
\begin{equation}
|\alpha|=2|g^{\rm T}(Q_{Y2}^{\rm T}B_mQ_{X2})g|\geq
2\sigma_{\min}(Q_{Y2}^{\rm T}B_mQ_{X2}),\label{alpha}
\end{equation}
where $\sigma_{\rm min}(C)$ denotes the smallest singular value of a
matrix $C$. Note that $P_mQ_Y$ and $Q_mQ_X$ form orthogonal bases of
the left subspace $span\{\hat U\}\oplus
span\{\hat{U}_{\perp}\}=span\{P_m\}$ and the right subspace
$span\{\hat V\}\oplus span\{\hat{V}_{\perp}\}=span\{Q_m\}$.
Exploiting (\ref{projmatrix}), (\ref{partition}), (\ref{perp}) and
(\ref{newbasis}) and keeping in mind that $span\{
P_mQ_{Y1}\}=\hat{U}$ and $span\{Q_mQ_{X1}\}=span\{\hat{V}\}$, we
obtain from (\ref{Aorth}) that the projection matrix of $A$ with
respect to $P_mQ_Y$ and $Q_mQ_X$ is
\begin{eqnarray*}
(P_mQ_Y)^{\rm T}A(Q_mQ_X)&=&\left(\begin{array}{c} Q_{Y1}^{\rm
T}P_m^{\rm T}\\
Q_{Y2}^{\rm T}P_m^{\rm T}
\end{array}
\right)A(Q_mQ_{X1}, Q_mQ_{X2})\\
&=&\left(\begin{array}{c} Q_{Y1}^{\rm
T}P_m^{\rm T}\\
\hat{U}_{\perp}^{\rm T}
\end{array}
\right)A(Q_mQ_{X1}, \hat{V}_{\perp})\\
&=&\left(\begin{array}{cc}
Q_{Y1}^{\rm T}B_mQ_{X1}&0 \\
0&\hat{U}_{\perp}^{\rm T}A\hat{V}_{\perp}
\end{array} \right)\\
&=&\left(\begin{array}{cc}
Q_{Y1}^{\rm T}B_mQ_{X1}&0 \\
0&Q_{Y2}^{\rm T}B_mQ_{X2}
\end{array}
\right),
\end{eqnarray*}
whose singular values $\tilde{\theta}_i,\,i=1,2,\ldots,m$, labeled
in increasing order, are the union of the singular values of
$Q_{Y1}^{\rm T}B_mQ_{X1}$ and $Q_{Y2}^{\rm T}B_mQ_{X2}$ and are just
the Ritz values of $A$ with respect to the left and right subspaces
$span\{P_m\}$ and $span\{Q_m\}$. By the singular value interlacing
property, we have $\sigma_i\leq\tilde{\theta}_i,\,i=1,2,\ldots,m$.
Furthermore, note that $Q_{Y1}^{\rm T}B_mQ_{X1}$ and $Q_{Y2}^{\rm
T}B_mQ_{X2}$ are the projection matrices of $A$ with respect to the
left subspaces $span\{\hat  U\}$ and $span\{\hat U_{\perp}\}$ and
the right subspaces $span\{\hat V\}$ and $span\{\hat V_{\perp}\}$,
respectively. Therefore, the singular values of $Q_{Y1}^{\rm
T}B_mQ_{X1}$ are $\tilde{\theta}_i,\,i=1,2,\ldots,k$ and approximate
the $k$ desired smallest singular values $\sigma_i$'s from above,
while the singular values of $Q_{Y2}^{\rm
T}B_mQ_{X2}=\hat{U}_{\perp}^{\rm T}A\hat{V}_{\perp}$ are
$\tilde{\theta}_i,\ i=k+1,\ldots,m$ and approximate
$\sigma_{k+1},\ldots,\sigma_{m}$ from above too. In particular, we
have
$$
\tilde{\theta}_{k+1}\geq\sigma_{k+1}.
$$
So it holds that
$$
|\alpha|\geq 2\sigma_{k+1}.
$$

Next we estimate $\beta$. We obtain from (\ref{gm}), (\ref{perp})
and (\ref{projmatrix})
\begin{eqnarray*}
\beta&=&g^{\rm T}(Q_{Y2}^{\rm T}B_mB_m^{\rm
T}Q_{Y2})g+\beta_m^2g^{\rm T}(Q_{Y2}^{\rm T}e_me_m^{\rm
T}Q_{Y2})g+g^{\rm T}(Q_{X2}^{\rm T}B_m^{\rm T}B_mQ_{X2})g\\
&=&g^{\rm T}(\hat{U}_{\perp}^{\rm T}AA^{\rm T}\hat{U}_{\perp})g+
g^{\rm T}(\hat{V}_{\perp}^{\rm T}A^{\rm T}A\hat{V}_{\perp})g\\
&\geq&\sigma_{\min}(\hat{U}_{\perp}^{\rm T}AA^{\rm
T}\hat{U}_{\perp})+\sigma_{\min}(\hat{V}_{\perp}^{\rm T}A^{\rm
T}A\hat{V}_{\perp}).
\end{eqnarray*}
Observe that $\hat{U}_{\perp}^{\rm T}AA^{\rm T}\hat{U}_{\perp}$ and
$\hat{V}_{\perp}^{\rm T}A^{\rm T}A\hat{V}_{\perp}$ are the
projection matrices of $AA^{\rm T}$ and $A^{\rm T}A$ with respect to
$span\{\hat{U}_{\perp}\}$ and $span\{\hat{V}_{\perp}\}$,
respectively. So their eigenvalues approximate some of the
eigenvalues $\sigma_{k+1}^2,\ldots,\sigma_N^2$ of $AA^{\rm T}$ and
$A^{\rm T}A$. Furthermore, since
$$
\hat{V}_{\perp}^{\rm T}A^{\rm
T}A\hat{V}_{\perp}-(\hat{U}_{\perp}^{\rm T}A\hat{V}_{\perp})^{\rm
T}(\hat{U}_{\perp}^{\rm T}A\hat{V}_{\perp})=\hat{V}_{\perp}^{\rm
T}A^{\rm T}(I-\hat{U}_{\perp}\hat{U}_{\perp}^{\rm
T})A\hat{V}_{\perp}
$$
is symmetric nonnegative definite, the smallest eigenvalue of $
\hat{V}_{\perp}^{\rm T}A^{\rm T}A\hat{V}_{\perp}$ is no less than
the smallest eigenvalue $\tilde{\theta}_{k+1}^2$ of
$(\hat{U}_{\perp}^{\rm T}A\hat{V}_{\perp})^{\rm
T}(\hat{U}_{\perp}^{\rm T}A\hat{V}_{\perp})$. As a consequence, it
follows from $\tilde{\theta}_{k+1}\geq\sigma_{k+1}$ that the
smallest eigenvalue of $\hat{V}_{\perp}^{\rm T}A^{\rm
T}A\hat{V}_{\perp}$ is bounded below by $\sigma_{k+1}^2$. Similarly,
we have
$$
\hat{U}_{\perp}^{\rm T}AA^{\rm
T}\hat{U}_{\perp}-(\hat{U}_{\perp}^{\rm
T}A\hat{V}_{\perp})(\hat{U}_{\perp}^{\rm T}A\hat{V}_{\perp})^{\rm
T}=\hat{U}_{\perp}^{\rm T}A(I-\hat{V}_{\perp}\hat{V}_{\perp}^{\rm
T})A^{\rm T}\hat{U}_{\perp},
$$
which is symmetric nonnegative definite. As a result, the smallest
eigenvalue of $\hat{U}_{\perp}^{\rm T}AA^{\rm T}\hat{U}_{\perp}$ is
no less than the smallest eigenvalue $\tilde{\theta}_{k+1}^2$ of
$(\hat{U}_{\perp}^{\rm T}A\hat{V}_{\perp})(\hat{U}_{\perp}^{\rm
T}A\hat{V}_{\perp})^{\rm T}$ and is bounded below by
$\sigma_{k+1}^2$ too. Therefore, we get
$$
\beta\geq 2\sigma_{k+1}^2.
$$
So it follows from (\ref{alphabound}) that the $\frac{1}{\xi}$'s
condition number $\nu$ in (\ref{condnum}) satisfies
$$
\nu\leq\frac{1}{2\sigma_{k+1}\sqrt{1+\sigma_{k+1}^2}},
$$
which is (\ref{nubound}).

If $\sigma_{k+1}$ is not very small, then $\nu$ is not large and all
$\xi$'s are not too small. Keep in mind the comments on
(\ref{backerror}) and (\ref{chi}). It is then clear that numerically
the proposed procedure can compute the refined harmonic shifts
$|\xi_i|$'s with relative accuracy $O(\epsilon_{\rm mach})$. \qquad
\endproof

We should point out that Theorem~\ref{theo10} holds without any
assumption on $span\{Q_m\}$ and $span\{P_m\}$, as is clearly seen
from the proof.

\subsection{Adaptive shifting strategy}

It has been observed \cite{larsen1} that if the $k$-th desired
$\sigma_k$ is very near a shift then IRLB with the exact shifts (the
unwanted Ritz values) converges very slowly and even stagnates. This
is also the case for IRRHLB with the refined harmonic shifts and
IRHLB with the harmonic shifts. The reason is that if some shift
$\mu_i$ is very near $\sigma_k$ then the new starting vector $q_1^+$
will nearly annihilate the component of the desired $v_k$, so that
the new subspace ${\cal K}_m(A^{\rm T}A,q_1^+)$ contains very little
information on $v_k$ and $\rho_k$ converges to $\sigma_k$ very
slowly or not at all.

In order to overcome this problem, for IRLB with the exact shifts,
Larsen \cite{larsen1} proposes an adaptive shifting strategy for
computing the largest singular triplets. He simply replaces a bad
shift to be defined below by a zero shift.  Jia and Niu
\cite{jianiu03} adapt it to IRRLB for computing the largest singular
triplets but modify it for computing the smallest singular triplets.
Their strategy works for IRHLB and IRRHLB: Define the relative gaps
of $\rho_k$ and all the shifts $\mu_i,i=1,2,\ldots,l$ by
  \begin{equation}
    {\rm  relgap}_{ki}=\left |\frac{(\rho_k-\varepsilon_k)-\mu_i}{\rho_k}\right |,
    \label{harmrelgap}
    \end{equation}
where $\varepsilon_k$ is the residual norm (\ref{harmstop}) or
(\ref{refinedstop}). We should note that $\rho_k-\varepsilon_k$ is
an approximation to $\sigma_k$. If ${\rm  relgap}_{ki}\leq 10^{-3}$,
$\mu_i$ is a bad shift and should be replaced by a suitable
quantity.

Expand $q_1$ as a linear combination of the right singular vectors
$\{v_j\}_{j=1}^N$:
$$
  q_1=\sum_{j=1}^N\alpha_jv_j.
$$
Then for the harmonic shifts $\mu_i=\theta_{k+i},\,i=1,2,\ldots,l$
we have from (\ref{update})
\begin{eqnarray*}
    \gamma q_1^+&=&\prod_{i=k+1}^m(A^{\rm T}A-\theta_i^2 I)q_1\\
    &=&\sum_{j=1}^k\alpha_j\prod_{i=k+1}^m(\sigma_j^2-\theta_i^2)v_j
    +\sum_{j=k+1}^N\alpha_j\prod_{i=k+1}^m(\sigma_j^2-\theta_i^2)v_j,
\end{eqnarray*}
So if $\theta_{k+1}$ is very near $\sigma_k$, which is the case that
$\sigma_{k+1}$ is very near $\sigma_k$, then $q_1^+$ has a very
small component in the direction of $v_k$. A good strategy is to
replace $\theta_{k+1}$ by the largest one among all the shifts, as
this strategy amplifies the components of $q_1^+$ in
$v_i,i=1,2,\ldots,k$ and meanwhile dampens those in
$v_i,i=k+1,\ldots,N$.

The above strategy applies to the refined harmonic shifts as well.

We now present IRHLB with the harmonic shifts and IRRHLB with the
refined harmonic shifts, respectively.
\smallskip

{\bf Algorithm 1. \ IRHLB with the harmonic shifts}

\begin{enumerate}
       \item Given a unit length starting vector $q_1$ of dimension $N$, the
       steps $m$, the number $k$ of the desired singular triplets
       and the convergence tolerance $tol$.

       \item Run the $m$-step Lanczos bidiagonalization process and construct
       $B_m,\,P_{m}$ and $Q_{m}$.

       \item Calculate the triplets
       $(\theta_i,\tilde s_i,\tilde w_i),\,i=1,2,\ldots,m$, by computing the
       singular values and right singular vectors of (\ref{svd1}) and by
       solving (\ref{bidiag}) and normalizing the solutions,
       and use the Rayleigh quotients $\rho_i=\tilde u_i^{\rm T}A\tilde v_i
       =\tilde s_i^{\rm T}B_m\tilde w_i$ as approximations to $\sigma_i,\,
       1=1,2,\ldots,k$.

       \item Replace $\theta_i$ by $\rho_i$ in (\ref{harmstop}).
       For $i=1,2,\ldots,k$, test if (\ref{harmstop}) is satisfied. If yes,
       compute $\tilde u_i$ and $\tilde v_i$ explicitly and stop.

       \item Implicitly restart the Lanczos bidiagonalization process using
       the harmonic shifts and the adaptive shifting strategy.
\end{enumerate}

{\bf Algorithm 2.\ IRRHLB with the refined harmonic shifts}

    \begin{enumerate}
       \item Given a unit length starting vector $q_1$ of dimension $N$, the
       steps $m$, the number $k$ of the desired singular triplets
       and the convergence tolerance $tol$.

       \item Run the $m$-step Lanczos bidiagonalization process
       and construct $B_m,\,P_{m}$ and $Q_{m}$.

       \item Calculate the triplets
       $(\theta_i,\tilde s_i,\tilde w_i),\,i=1,2,\ldots,m$ by computing
       the singular values and right singular vectors of (\ref{svd1})
       and by solving (\ref{bidiag}) and normalizing the solutions, and
       use the Rayleigh quotients $\rho_i
       =\tilde u_i^{\rm T}A\tilde v_i=\tilde s_i^{\rm T}B_m\tilde
       w_i$ as approximations to $\sigma_i,\,1=1,2,\ldots,k$.

       \item For each $\rho_i,i=1,2,\ldots,k$, compute $\hat x_i$ and $\hat y_i$
       in Theorem~\ref{theo4}.

       \item For $i=1,2,\ldots,k$, test if (\ref{refinedstop}) is satisfied.
       If yes, compute $\hat u_i$ and $\hat v_i$ by (\ref{refinevector})
       explicitly and stop.

       \item Implicitly restart the Lanczos bidiagonalization process
       using the refined harmonic shifts and the adaptive shifting strategy.
       \end{enumerate}

\section {Numerical experiments}

We have developed the experimental Matlab codes of IRRHLB, IRHLB,
IRRLB and IRLB. The latter two were named IRRBL and IRBL in
\cite{jianiu03} and were originally developed based on the lower
Lanczos bidigonalization process. Here we have developed their upper
Lanczos bidiagonalization versions. These four codes call the upper
Lanczos bidiagonalization process in Baglama and Reichel's code
IRLBA, and some parameters and defaults are the same as those used
in IRLBA. We compare IRRHLB with IRRHLB, IRRLB, IRLB, IRLBA and
IRLANB in this section and report numerical results. Numerical
experiments were run on an Intel Core 2 E6320 with CPU 1.86GHz and
RAM 2GB under the Window XP operating system using Matlab 7.1 with
$\epsilon_{\rm mach}=2.22\times 10^{-16}$. The stopping criteria are
$$
   stopcrit=\max_{1 \leq i \leq k}\sqrt{\|A\tilde v_i-\rho_i
    \tilde u_i\|^2+\|A^{\rm T}\tilde u_i-\rho_i\tilde v_i\|^2} \ \ (\rm  IRHLB)
$$
and
$$
    stopcrit=\max_{1 \leq i \leq k}\sqrt{\|A\hat v_i-\rho_i
    \hat u_i\|^2+\|A^{\rm T}\hat u_i-\rho_i\hat v_i\|^2}\ \ (\rm
    IRRHLB).
$$
If
  \begin{equation}
    stopcrit=\frac{stopcrit} {\|A\|}<tol, \label{stop}
  \end{equation}
then stop. Similar criteria apply to IRLB and IRRLB as well. In
(\ref{stop}), $\|A\|$ is replaced by the maximum of the current
largest (harmonic) Ritz value and the old one obtained at last
restart. Some parameters in IRRHLB, IRHLB, IRRLB and IRLB are
described in Table~\ref{description}.

\begin{table} [htp]
\caption {Parameters in IRRHLB, IRRLB, IRHLB and IRLB}
\label{description}
\begin{center}
    \begin{tabular}{|ll|} \hline
        Parameters & Description\\
        $k$      & Number of the desired singular triplets. \\
                 & Default value: $k=6$.\\
        $adjust$ & Number added to $k$ to speed up convergence.\\
                 & Default value: $adjust=3$. \\
        $disps$  & When $disps>0$, the $k$ desired approximate singular values\\
                 & and norms of associated absolute residual error are displayed\\
                 & each iteration. $disps=0$ inhibits display of these quantities.\\
                 & Default value: $disps=1$.\\
        $M\_B$   & Maximum of Lanczos bidiagonalization steps.  \\
                 & Default value: $M\_B=20$.\\
        $maxit$  & Maximum number of restarts.\\
                 & Default value: $maxit=300$.\\
        $sigma$  & A 2-letter string which specifies which extreme singular
                   triplets are to\\
                 & be computed, 'SS' for the smallest and 'LS' for the largest.\\
                 & Default value: $sigma=$'SS'.\\
        $tol$    & User defined relative tolerance to check convergence.\\
                 & Default value: $tol=10^{-6}$.\\
        $v_0$    & min$(M,N)$-dimensional initial vector of
                  Lanczos bidiagonalization.\\
                 & Default value: $v_0=randn(\min(M,N),1)$. \\\hline
    \end{tabular}
    \end{center}
\end{table}

For large matrix eigenproblems, in order to speed up convergence,
ARPACK ({\sf eigs}) and Implicitly Restarted Refined Arnoldi
Method (IRRA) \cite{jia99a} compute $k+3$ approximate eigenpairs,
so the number of shifts is $m-(k+3)$ when $k$ eigenpairs are
desired. This strategy adapts to Krylov type subspace algorithms
for SVD problems, for instance, the default parameter $adjust=3$
in IRLBA, which means that the updating subspaces are augmented
with $k+3$ Ritz or harmonic Ritz vectors. With $m-(k+3)$ exact
shifts and harmonic shifts used, IRLB, IRHLB and IRLANB retain
$k+3$ Ritz vectors and harmonic Ritz vectors in the updating
subspaces, respectively.

We mention that our codes IRRHLB, IRHLB, IRRLB and IRLB as well as
IRLBA and IRLANB (we used the newest available code ${\sf
irlanb_{-}review}$) do not involve any shift-and-invert matrix when
computing the smallest singular triplets, while the Matlab internal
function {\sf svds} needs to factorize $\tilde A$ of
(\ref{argumented}). In this context, we assume that $A$ is too large
to allow any factorization of $\tilde A$ due to excess memory and/or
computational cost, so we do not compare the above six algorithms
with {\sf svds}.

Our experiments consist of three subsections. In the first two
subsections, we test IRRHLB on a set of matrices having the
clustered smallest singular values and on a set of ill-conditioned
matrices, respectively. We show that IRRHLB works well on them and
confirm some theory. In the third subsection, we compare IRRHLB with
the five other algorithms on seven practical problems that include
very difficult, difficult and general ones, illustrating that IRRHLB
is at least competitive with and can be much more efficient than the
five other ones.

\subsection{IRRHLB for the clustered smallest singular values}

\begin{figure}[t]
\epsfig{figure=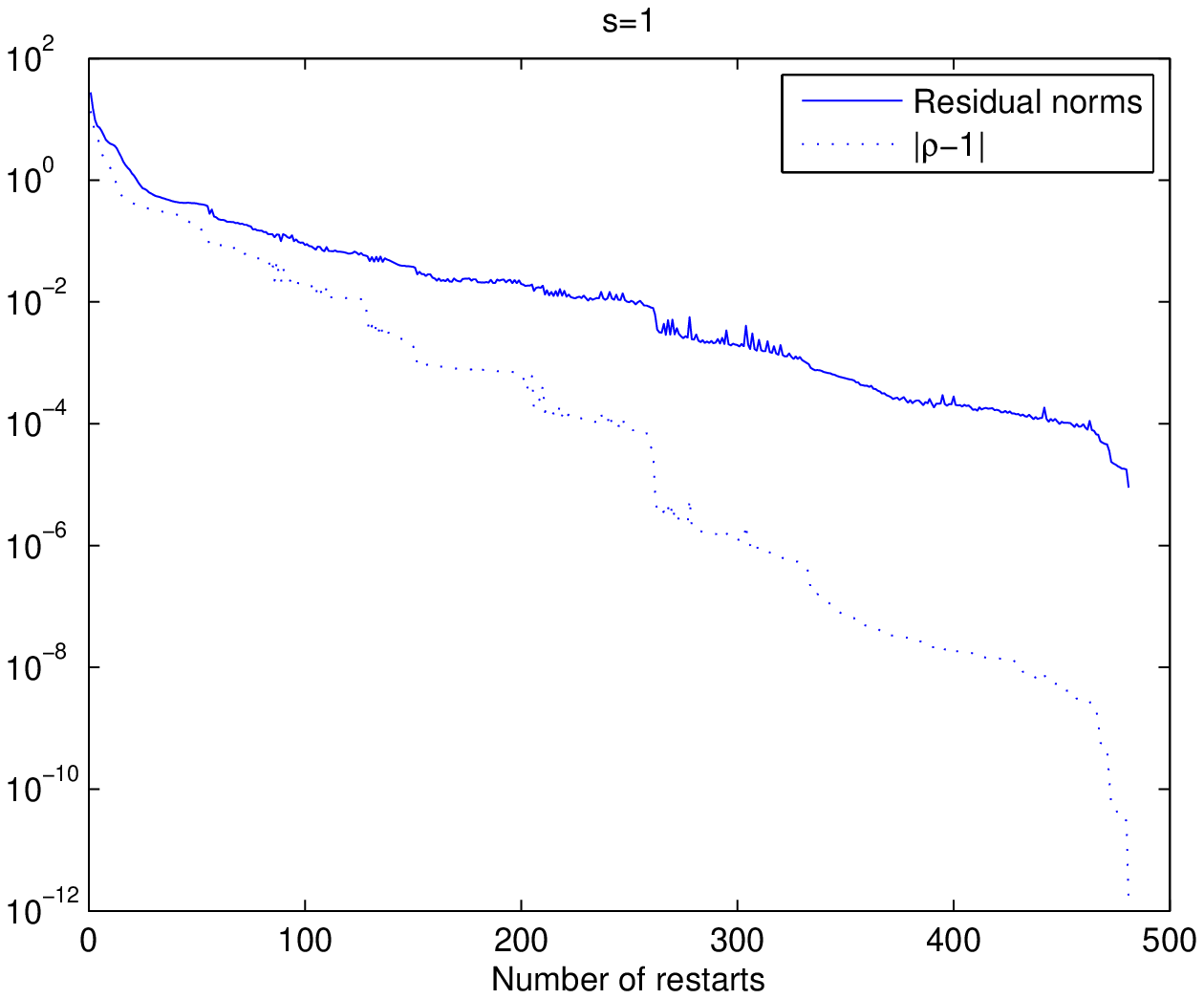,height=2.0in}
\epsfig{figure=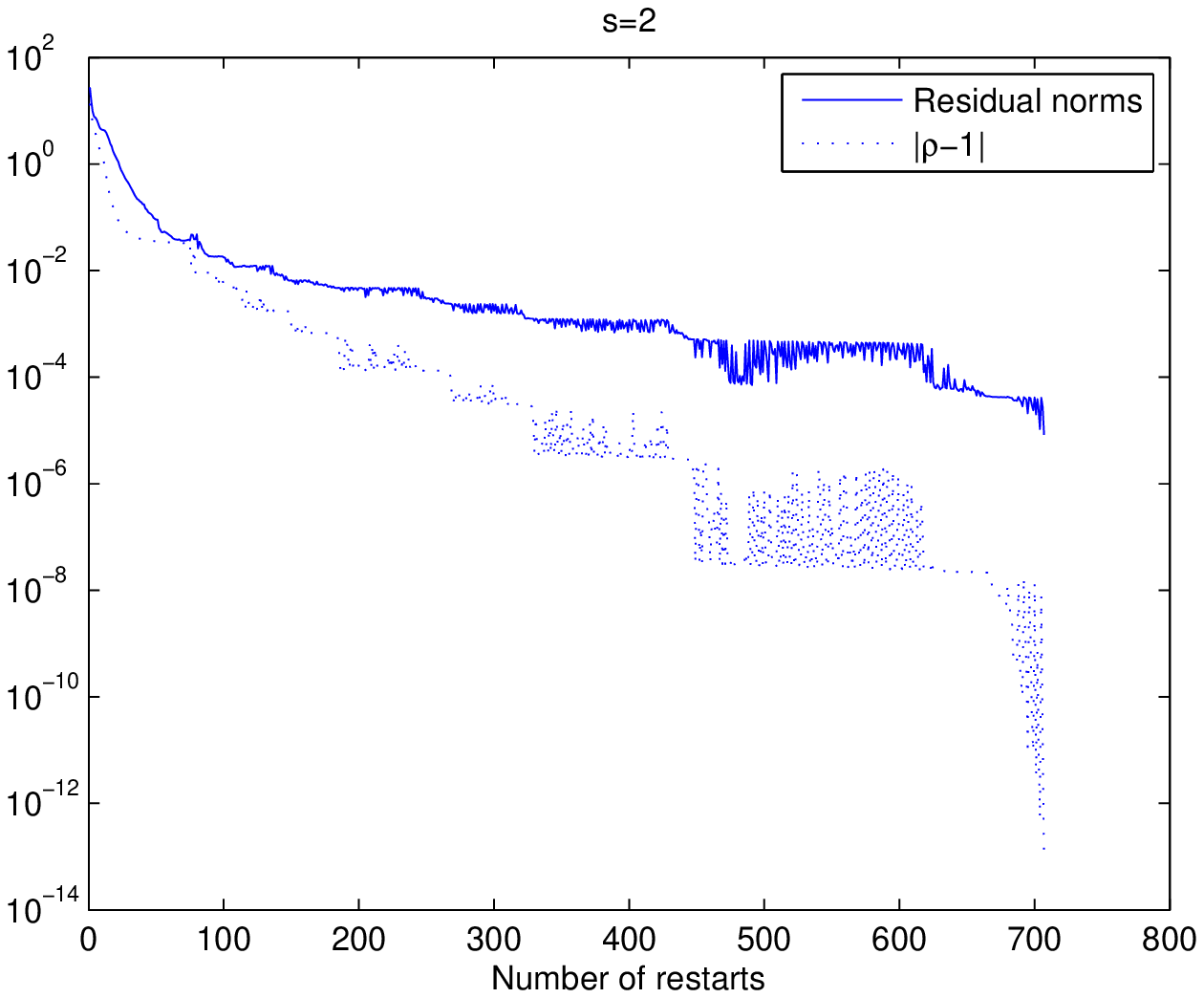,height=2.0in}
\epsfig{figure=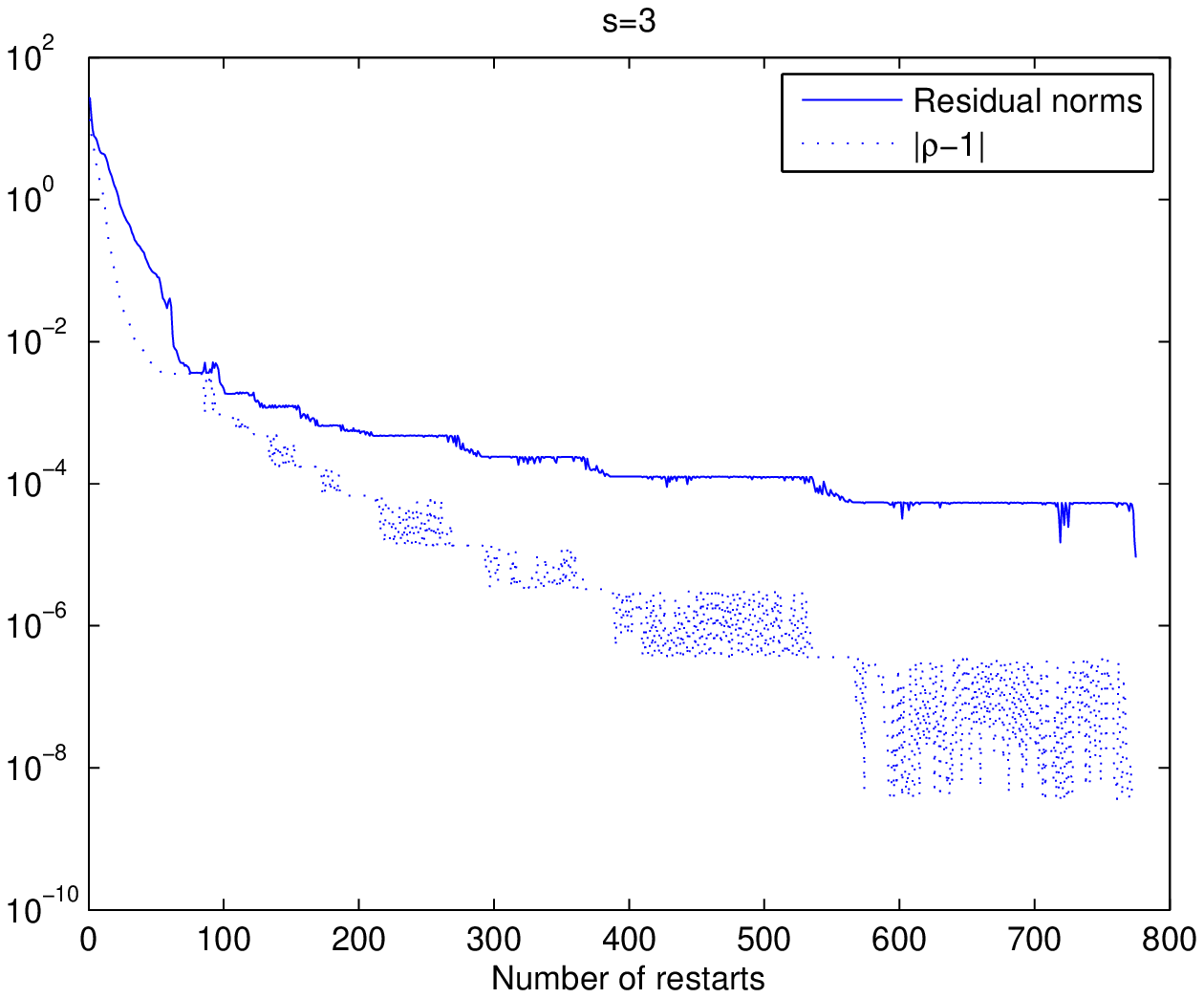,height=2.0in}
\epsfig{figure=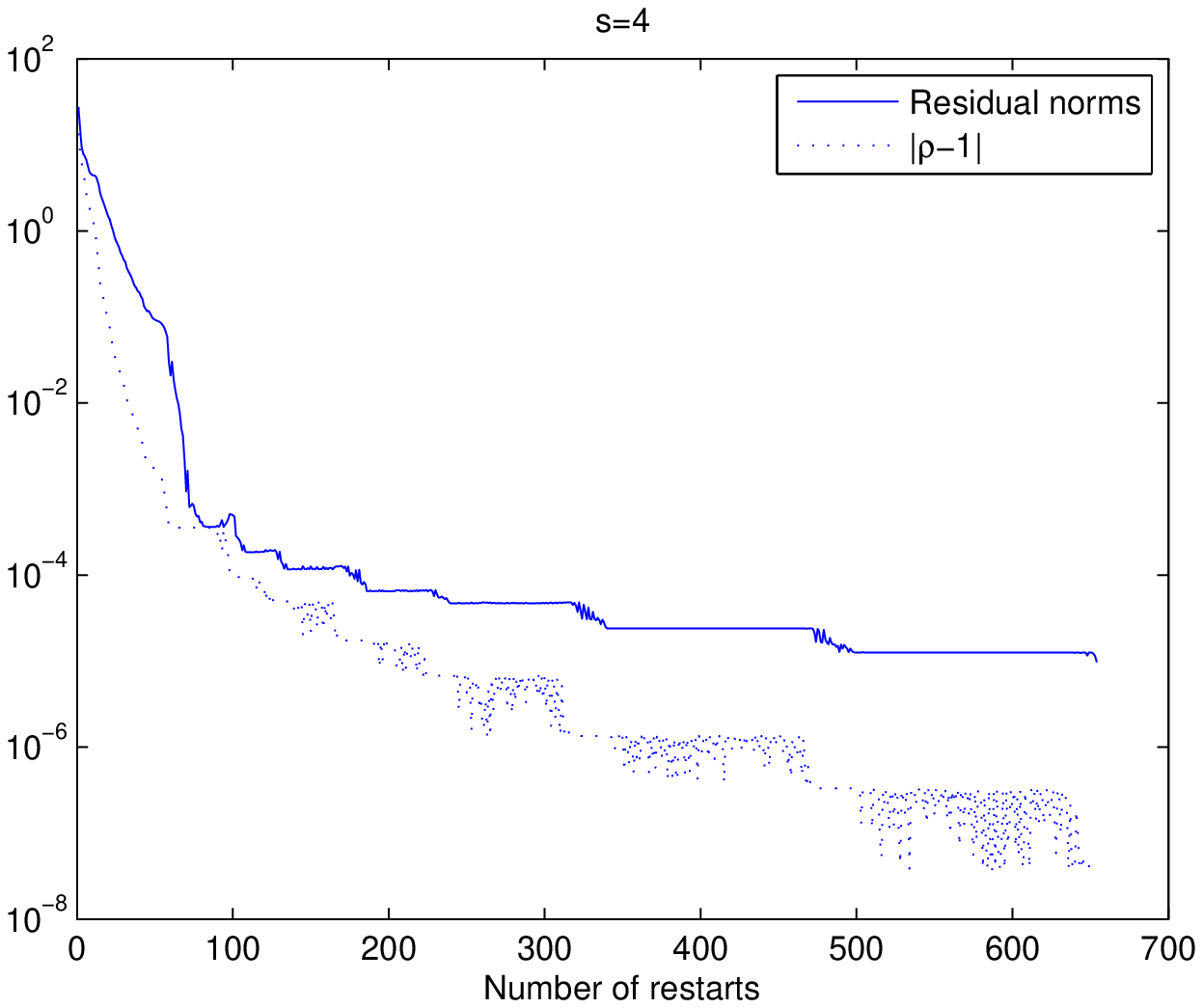,height=2.0in}
\caption{Experiments with
the matrices (\ref{Ex81}), with the clustering smallest singular
value. Solid lines correspond to residual norms. Dotted lines to
relative errors of the smallest singular value computed by
IRRHLB.}\label{fig1}
\end{figure}

This set of experiments is designed to see how IRRHLB behaves for
the clustered smallest singular values. Similar to those matrices in
\cite{Kokio04}, we constructed a sequence of diagonal matrices
$A_s\in {\cal R}^{n\times n}$, $n=1000,\ s=1,\ldots,4$ whose nine
smallest singular values become increasingly more clustered as $s$
increases. In the Matlab language: {\small
\begin{equation}\label{Ex81}
A_s=spdiags([1:10^{-s}:1+9*10^{-s},2:1:1000]',0,1000,1000],s=1,\ldots,4,
\end{equation}}
whose smallest singular value $\sigma_1=1$ and $\kappa(A_s)=991$ for
all $s$. Since $\kappa(A_s)$ is moderate, it is expected that IRRHLB
computes $\sigma_1$ accurately if it works. We computed $\sigma_1$
by taking the parameters
$$
opts.m=50, opts.maxit = 2000, opts.adjust=9, opts.tol=1e-8
$$
and using the same starting vector generated randomly in a normal
distribution for all $s$. Figure~\ref{fig1} plots absolute residual
norms of the computed singular triplets and relative errors
$|\rho-1|$, respectively. We see that IRRHLB succeeded for all $s$
and computed the smallest singular value accurately. In the worst
case $s=4$, IRRHLB gave the relative error $|\rho-1|=3.7\times
10^{-8}$, the same order as the backward error that equals the
relative residual norm $10^{-8}$. For $s=3$, the relative error
$|\rho-1|=3.5\times 10^{-9}$. Both are in agreement with the
standard perturbation theory \cite{golub96,stewart}. For $s=1,2$,
the relative errors $|\rho-1|$ are $1.3\times 10^{-12}$ and
$8.1\times 10^{-14}$, respectively, a few order smaller than the
predicted relative error $O(10^{-8})$. Also, IRRHLB used
considerably fewer restarts for $s=1$ than for the other bigger $s$
but had comparable restarts for $s=2,3,4$. It was expected that
IRRHLB converged faster for $s=1$ than for $s=2,3,4$ since the gap
of $\sigma_1$ and $\sigma_2$ is considerably bigger for $s=1$ than
those for $s=2,3,4$.

\subsection{IRRHLB for ill-conditioned matrices}

\begin{figure}[htp]
\epsfig{figure=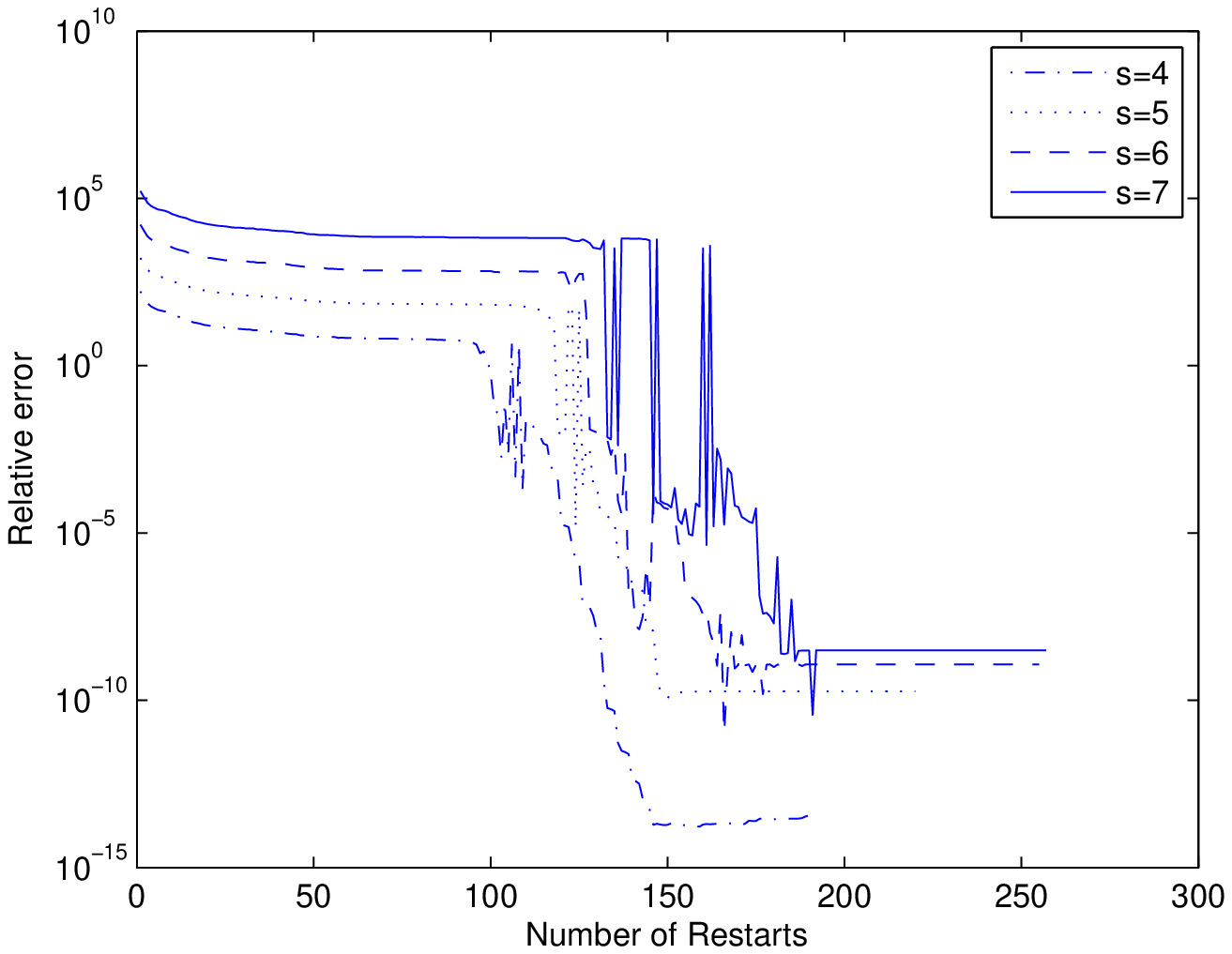,height=2.0in}
\epsfig{figure=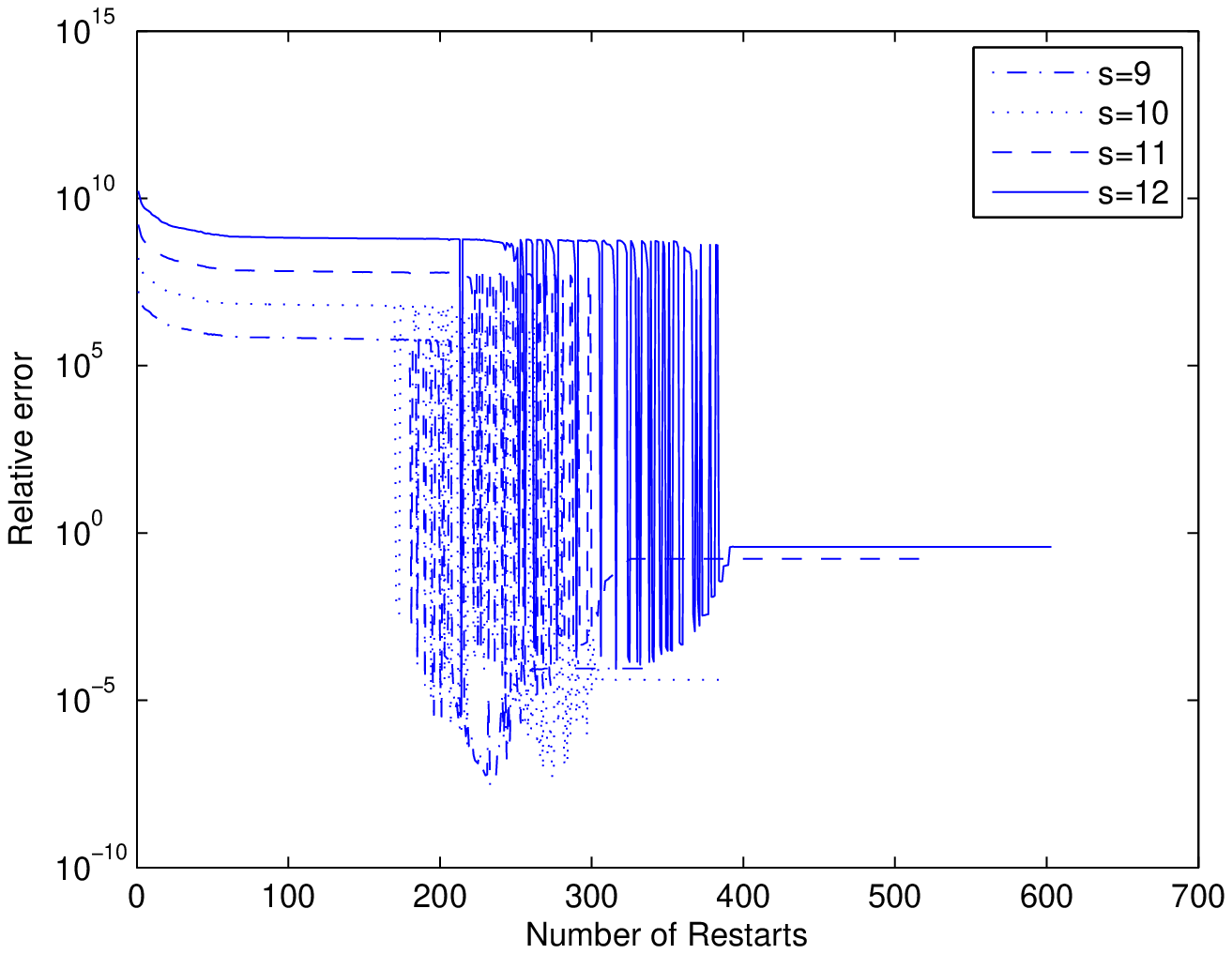,height=2.0in} \caption{Experiments
with the increasingly ill-conditioned diagonal matrices
(\ref{Ex82}),with
$\kappa(A_s)=10^s,s=4,5,6,7,9,10,11,12$.}\label{fig2}
\end{figure}

We investigate the behavior of IRRHLB for a set of ill-conditioned
matrices. Similar to those matrices in \cite{Kokio04}, we
constructed a sequence of bidiagonal matrices $A_s\in {\cal
R}^{n\times n}$, $n=1000,\,s=4,\ldots,7$ and $9,\ldots,12$ with
increasing condition numbers:
\begin{equation}\label{Ex82}
A_s=spdiags(linspace(1,10^s,1000)',0,1000,1000),
\end{equation}
whose smallest singular value $\sigma_1=1$ and condition numbers
$\kappa(A_s)=10^s$. We computed $\sigma_1$ by taking the parameters
$$
opts.m=50, opts.maxit = 2000, opts.adjust=3, opts.tol=1e-14
$$
using the same starting vector generated randomly in a normal
distribution for all $s$. Figure~\ref{fig2} plots relative errors
$|\rho-1|$. It was seen from the figure that IRRHLB computed the
smallest singular value with relative error smaller than $10^{-8}$
for $s=4,5,6,7$. This confirms the perturbation theory: the smaller
$s$ is, the smaller the relative error is. For more ill-conditioned
cases $s=9,10,11,12$, the accuracy of the computed smallest singular
values deteriorated significantly. For $s=9,10$, the relative errors
are $9.0\times 10^{-5}$ and $4.0\times 10^{-5}$, and the computed
singular values have four and five correct decimal digits,
respectively. In the worst-conditioned case that $s=12$, the
relative error is $0.3848$; for $s=11$, the relative error is
$0.1678$, and the computed smallest singular value was a little bit
more accurate than that for $s=12$.

We also tested $opts.tol=1e-12$ and compared the results with those
for $opts.tol=1e-14$. We found that although residual norms
continued decreasing until $10^{-14}$, the accuracy of the computed
singular values was not improved further as residual norms decreased
from $opts.tol=1e-12$ to $opts.tol=1e-14$. This was reflected by the
figures, where we saw that relative errors did not decrease further
and stabilized, starting from some restart for each $s$ except
$s=9$. The curves for $s=9$ jumped up and down when the algorithm
approached convergence, but kept below $10^{-4}$. All these are in
accordance with the predicted relative errors, which should not be
bigger than a very modest multiple of $\kappa(A_s)\times opts.tol$.
Another important observation is that IRRHLB used more restarts as
$s$ increases. Since the ratio
$\frac{\sigma_N-\sigma_1}{\sigma_2-\sigma_1}$, the spread over the
gap of $\sigma_1$ and $\sigma_2$, increases as $s$ does, it is more
difficult for IRRHLB to solve the SVD problem as $s$ increases. We
also saw that the curves of relative errors oscillated quite often
in the middle of convergence processes. A careful observation
revealed that IRRHLB started to be on its way to compute the desired
singular value at some stage but lost it soon. Then it adaptively
adjusted convergence repeatedly and eventually was on the correct
way to converge. These phenomena may be explained by
Theorem~\ref{theo1} and the comments followed.

\subsection{Experiments of the six algorithms on practical problems}

We now do numerical experiments on several selected problems that
include very difficult, difficult and general ones. We compare
IRRHLB with the five other algorithms: IRRLB, IRHLB, IRLB, IRLBA
and IRLANB.

Table~\ref{table1} lists seven test matrices from \cite{bai,duff}
and some of their basic properties. Except well1852, all other
matrices are square matrices. Note that the ratio $\frac{{\rm
spread}(A)}{{\rm gap}(k)}$ indicates whether or not the six
algorithms are difficult to converge. The bigger $\frac{{\rm
spread}(A)}{{\rm gap}(k)}$ is, the more slowly the algorithms
converge generally. From the table, we see that the $k(=1,3,5,10)$
desired smallest singular values of all the matrices are quite
clustered; among them the matrices fidap4, jagmaesh8 and lshp3205
are the most difficult, the matrix plat1919 is relatively
difficult, and the matrices well1850 and dw2048 are general. We
see that all $\kappa(A)$'s are not very large, so the columns of
$A$ are strongly linearly independent. It is expected that if the
algorithms converge then they can compute the smallest singular
values with relative errors no more than a very modest multiple of
$\kappa(A)\times opts.tol$.

\begin{table}[htb]
\caption {Six test matrices: fidap4 of order $1601\times 1601$,
jagmesh8 of order $1141\times 1141$, lshp3205 of order $3205
\times 3205$, well1850 of order $1850\times 712$, pde2961 of order
$2961\times 2961$ and plat1919 of order $1919\times 1919$. ${\rm
spread}(A)=\sigma_N-\sigma_1$ and ${\rm gap}(k)=\min
(\sigma_{i+1}-\sigma_i), i=1,2,\ldots,k.$ } \label{table1}
\begin{center}
{\small
\begin{tabular}{|c|c|c|c|c|c|c|c|} \hline
Matrix&fidap4&jagmesh8&lshp3025&well1850&dw2048&pde2961&plat1919\\\hline
nnz($A$)&31837&7465&20833&8755&10114&14580&32399\\ \hline
$\kappa(A)$&$2.4e+3$&$2.2e+4$&$6.8e+4$&$1.1e+2$&$2.1e+3$&$6.4e+2$&$3.7e+2$\\\hline
spread($A$)&$1.6e+0$&$6.8e+0$&$7.0e+0$&$1.8e+0$&$9.8e-1$&$1.0e+1$&$2.3e+0$\\\hline
gap(1)&$1.5e-3$&$1.7e-3$&$1.8e-3$&$3.0e-3$&$2.6e-3$&$8.2e-3$&$2.6e-3$\\\hline
gap(3)&$2.5e-4$&$1.6e-3$&$9.1e-4$&$3.0e-3$&$2.9e-4$&$2.4e-3$&$1.8e-3$\\\hline
gap(5)&$2.5e-4$&$4.8e-5$&$1.8e-4$&$3.0e-3$&$2.9e-4$&$2.4e-3$&$2.7e-4$\\\hline
gap(10)&$2.5e-4$&$4.8e-5$&$2.2e-5$&$2.6e-3$&$1.6e-4$&$5.2e-4$&$2.0e-4$\\\hline
\end{tabular}}
\end{center}
\end{table}

We computed the $k$ smallest singular triplets for different $k$.
To make a reasonable comparison, for each matrix except well1850
we used the same starting vector generated randomly in a normal
distribution for the six algorithms. For well1850, we took the
same starting vector $u_0={\sf randn}(1850,1)$ in IRLANB and the
same starting vector $v_0={\sf randn}(712,1)$ in the five other
algorithms. In all tables, denote by $iter$ the number of
restarts, by $time$ CPU time in second, by $n.c$ non-convergence
after 2000 restarts are used, and by $mv$ the number of
matrix-vector products. Since matrix-vector products involving $A$
are equal to those involving $A^{\rm T}$, we only count the number
of matrix-vector products involving $A$. We compare restarts and
matrix-vector products as well as CPU time needed by all the codes
for the same $k$ and $m$. The former two quantities reflect the
overall efficiency of the codes more fairly and reasonably.

By the above description, in IRHLB, IRRLB, IRHLB, IRLB and IRLBA
we took the input parameters
$$
opts.k = k,\ opts.M\_B = m,\ opts.tol = tol, \ opts.maxit = 2000, \
opts.v0 = v_0
$$
and the others as defaults. In IRLANB we took
\begin{eqnarray*}
&eignum = k, \ options.k = m-1,\ options.l = m-1-(k+3),&\\
&options.u0=u_0, \ options.maxit = 2000,\
options.version='harmonic'&
\end{eqnarray*}
 and the others as defaults. This
parameters make all the codes compute the approximate singular
triplets with respect to certain subspaces of the same dimension
$m$ and use the same number of shifts at each restart.

We found that fidap4, jagmesh8 and lshp3025 challenged most of the
six algorithms. Tables~\ref{tablefida}--\ref{tablelshp} report the
results obtained by IRRHLB and IRRLB for $opts.tol=1e-6$.

\begin{center}
\begin{table}[htp]
\caption {fidap4 for $k=1, 3, 5, 10$}\label{tablefida}
\begin{center}
    \begin{tabular}{|c|c|c|c|c|c|c|c|c|c|} \hline
      $k=1$& \multicolumn{3} {|c|} {$m=15$} & \multicolumn{3} {|c|}
      {$m=20$}& \multicolumn{3} {|c|} {$m=25$}\\  \hline
       Method & $iter$ & $time$ & $mv$ & $iter$ & $time$ & $mv$ & $iter$
      & $time$ & $mv$\\ \hline
      IRRHLB&807&46.5&8881&430&4428&6884&375&62.0&7879    \\ \hline
      IRRLB&$n.c$&-&-&$1929$&207&30868&1208&187&25372 \\ \hline
      $k=3$& \multicolumn{3} {|c|} {$m=15$} & \multicolumn{3} {|c|}
      {$m=20$}& \multicolumn{3} {|c|} {$m=25$}\\  \hline
       Method & $iter$ & $time$ & $mv$ & $iter$ & $time$ & $mv$ & $iter$
      & $time$ & $mv$\\ \hline
      IRRHLB&1083&55.4&9753&718&60.5&10058&447&64.1&8499    \\ \hline
      IRRLB &$n.c$&-&-&$n.c$&-&-&1526&213&29000 \\ \hline
      $k=5$& \multicolumn{3} {|c|} {$m=20$} & \multicolumn{3} {|c|}
      {$m=25$}& \multicolumn{3} {|c|} {$m=30$}\\  \hline
       Method & $iter$ & $time$ & $mv$ & $iter$ & $time$ & $mv$ & $iter$
      & $time$ & $mv$\\ \hline
      IRRHLB&1797&71.7&12587&1151&114&13820&777&100&13217   \\ \hline
      $k=10$& \multicolumn{3} {|c|} {$m=20$} & \multicolumn{3} {|c|}
      {$m=25$}& \multicolumn{3} {|c|} {$m=30$}\\  \hline
       Method & $iter$ & $time$ & $mv$ & $iter$ & $time$ & $mv$ & $iter$
      & $time$ & $mv$\\ \hline
      IRRHLB&1959&139&13726&1084&142&13021&737&134&12542    \\ \hline
    \end{tabular}
    \end{center}
\end{table}
\end{center}

\begin{center}
\begin{table}[htp]
\caption {jagmesh8 for $k=1, 3, 5, 10$}\label{tablejag}
\begin{center}
    \begin{tabular}{|c|c|c|c|c|c|c|c|c|c|} \hline
      $k=1$& \multicolumn{3} {|c|} {$m=20$} & \multicolumn{3} {|c|}
      {$m=25$} &\multicolumn{3} {|c|} {$m=30$}\\  \hline
       Method & $iter$ & $time$ & $mv$ & $iter$ & $time$ & $mv$ & $iter$ & $time$ & $mv$\\ \hline
      IRRHLB&$1167$&81.3&18676&953&122&20017&828&117&21532    \\ \hline
      $k=3$& \multicolumn{3} {|c|} {$m=20$} & \multicolumn{3} {|c|}
      {$m=25$} &\multicolumn{3} {|c|} {$m=30$}\\  \hline
       Method & $iter$ & $time$ & $mv$ & $iter$ & $time$ & $mv$& $iter$ & $time$ & $mv$ \\ \hline
      IRRHLB&1563&101&21888&1239&136&23547&897&126&21534    \\ \hline
      $k=5$& \multicolumn{3} {|c|} {$m=20$} & \multicolumn{3} {|c|}
      {$m=25$} &\multicolumn{3} {|c|} {$m=30$}\\  \hline
       Method & $iter$ & $time$ & $mv$ & $iter$ & $time$ & $mv$& $iter$ & $time$ & $mv$ \\ \hline
      IRRHLB&1521&89.6&18260&1108&135&11844&761&107&16750   \\ \hline
      $k=10$& \multicolumn{3} {|c|} {$m=25$} & \multicolumn{3} {|c|}
      {$m=30$} &\multicolumn{3} {|c|} {$m=35$}\\  \hline
       Method & $iter$ & $time$ & $mv$ & $iter$ & $time$ & $mv$& $iter$ & $time$ & $mv$ \\ \hline
      IRRHLB&1216&125&14605&882&148&15007&793&169&17459    \\ \hline
      IRRLB &$n.c$&-&-&1763&316&29984&810&175&17833 \\ \hline
      IRHLB &$n.c$&-&-&$n.c$&-&-&1707&204&37567 \\ \hline
      IRLB  &$n.c$&-&-&$n.c$&-&-&1853&220&40779 \\ \hline
      IRLBA &$n.c$&-&-&$n.c$&-&-&1919&37&42230 \\ \hline
    \end{tabular}
    \end{center}
\end{table}
\end{center}

\begin{center}
\begin{table}[ht]
\caption{lshp3025 for $k=1, 3, 5$}\label{tablelshp}
\begin{center}
    \begin{tabular}{|c|c|c|c|c|c|c|c|c|c|} \hline
      $k=1$& \multicolumn{3} {|c|} {$m=30$} & \multicolumn{3} {|c|}
      {$m=40$}& \multicolumn{3} {|c|} {$m=50$}\\  \hline
       Method & $iter$ & $time$ & $mv$ & $iter$ & $time$ & $mv$& $iter$ & $time$ & $mv$ \\ \hline
      IRRHLB&1116&293&30320&886&405&31900&947&918&43566    \\ \hline
      IRRLB &$n.c$&-&-&$n.c$&-&-&1604&1522&73788 \\ \hline
      $k=3$& \multicolumn{3} {|c|} {$m=30$} & \multicolumn{3} {|c|}
      {$m=40$}& \multicolumn{3} {|c|} {$m=50$}\\  \hline
       Method & $iter$ & $time$ & $mv$ & $iter$ & $time$ & $mv$ & $iter$
       & $time$ & $mv$ \\ \hline
      IRRHLB&1520&496&36486&1139&761&38732&971&978&42730    \\ \hline
      IRRLB &$n.c$&-&-&$n.c$&-&-&1116&11133&49110 \\ \hline
      $k=5$& \multicolumn{3} {|c|} {$m=40$} & \multicolumn{3} {|c|}
      {$m=50$}& \multicolumn{3} {|c|} {$m=60$}\\  \hline
       Method & $iter$ & $time$ & $mv$ & $iter$ & $time$ & $mv$ & $iter$
       & $time$ & $mv$ \\ \hline
      IRRHLB&$n.c$&-&-&1931&1906&81110&1656&2439&86120   \\ \hline
    \end{tabular}
    \end{center}
\end{table}
\end{center}

Clearly, for fidap4 and lshp3205, IRRHLB worked well and solved
the problem successfully while IRRLB only performed well in some
cases and was less efficient than IRRHLB. In contrast, IRHLB,
IRLB, IRLBA and IRLANB performed more poorly and they all failed
to converge for fidap4 and lshp3205. For jagmesh8, IRRHLB still
worked robustly and efficiently, but IRRLB succeeded only in a few
cases and IRHLB, IRLB, IRLBA and IRLANB behaved more badly. They
all were considerably less efficient than IRRHLB if they worked.
We found that, generally, the bigger $k$ was, the more restarts
IRRHLB and IRRLB used for the same $m$. This should not be
surprising as the problem for a bigger $k$ is generally more
difficult to solve than that for a smaller $k$. We also observed
that all the smallest singular values were computed with relative
errors no more than a modest multiple of $\kappa(A)\times
10^{-6}$.

We had more observations on the behavior of IRHLB, IRLB, IRLBA and
IRLANB on these three difficult problems. For example, the
residual norms obtained by them may oscillated but decreased very
slowly; they may have first decreased to some stage and then
oscillated; they might have first decreased, then stabilized and
did not decrease further; they might have decreased to some stage
and then increased. Therefore, IRRHLB is not only the best but
also the unique choice for fidap4, jagmesh8 and lshp3025 for most
of the given $k$'s and $m$'s.

We tested well1850, pde2961, dw2048 and plat1919 for
$opts.tol=1e-6,\,1e-9$, respectively.
Tables~\ref{tablewell6}--\ref{tableplat6} report the results for
$opts.tol=1e-6$. We do not list the corresponding results for
$opts.tol=1e-9$, as will be explained shortly.

\begin{table}[htp]
\caption {well1850 for $k=1, 3, 5, 10,
tol=1e-6$}\label{tablewell6}
\begin{center}
    \begin{tabular}{|c|c|c|c|c|c|c|c|c|c|} \hline
      \multicolumn{10}{|c|}{$k=1$}\\ \hline
      $m$&  \multicolumn{3} {|c|} {15}& \multicolumn{3} {|c|} {20}
      & \multicolumn{3} {|c|} {25}\\  \hline
      Algorithms & $iter$ & $time$ & $mv$ & $iter$ & $time$
      & $mv$& $iter$ & $time$ & $mv$\\ \hline
      IRRHLB &71&2.67&785&43&2.98&692&35&4.18&739   \\ \hline
      IRRLB  &69&2.51&763&62&4.48&996&35&4.16&739  \\ \hline
      IRHLB  &168&5.32&1852&83&4.99&1332&51&5.55&1075 \\ \hline
      IRLB   &183&5.99&2017&91&5.63&1460&55&6.01&1159  \\ \hline
      IRLBA  &191&1.90&2105&93&1.21&1492&57&1.00&1201 \\ \hline
      IRLANB &279&7.88&2795&133&6.61&2000&82&6.40&1645 \\ \hline
      \multicolumn{10}{|c|}{$k=3$}\\ \hline
      $m$&  \multicolumn{3} {|c|} {15}& \multicolumn{3} {|c|} {20}
      & \multicolumn{3} {|c|} {25}\\  \hline
      Algorithms & $iter$ & $time$ & $mv$ & $iter$ & $time$
      & $mv$& $iter$ & $time$ & $mv$\\ \hline
      IRRHLB &80&3.10&726&51&3.53&720&35&4.08&671   \\ \hline
      IRRLB  &103&4.15&933&63&4.43&888&41&5.16&785  \\ \hline
      IRHLB  &171&5.85&1545&76&4.51&1070&46&4.47&880 \\ \hline
      IRLB   &184&4.96&1662&82&4.84&1154&49&4.80&937  \\ \hline
      IRLBA  &189&1.70&1707&83&1.01&1168&50&0.86&956 \\ \hline
      IRLANB &259&6.11&2079&109&4.79&1424&63&4.10&1141 \\ \hline
      \multicolumn{10}{|c|}{$k=5$}\\ \hline
      $m$& \multicolumn{3} {|c|} {15}& \multicolumn{3} {|c|} {20}
      & \multicolumn{3} {|c|} {25}\\  \hline
      Algorithms & $iter$ & $time$ & $mv$ & $iter$ & $time$
      & $mv$& $iter$ & $time$ & $mv$\\ \hline
      IRRHLB &105&3.18&743&60&4.00&728&43&4.99&739   \\ \hline
      IRRLB  &161&5.53&1135&70&4.51&848&50&4.64&688  \\ \hline
      IRHLB  &248&6.43&1744&94&5.14&1136&53&4.06&909 \\ \hline
      IRLB   &275&6.08&1933&103&4.70&1244&58&4.69&994  \\ \hline
      IRLBA  &292&2.23&2052&108&1.26&1304&60&0.97&1028 \\ \hline
      IRLANB &388&7.30&2337&128&4.58&1417&69&3.41&1113 \\ \hline
      \multicolumn{10}{|c|}{$k=10$}\\ \hline
      $m$& \multicolumn{3} {|c|} {20}& \multicolumn{3} {|c|} {25}
      & \multicolumn{3} {|c|} {30}\\  \hline
      Algorithms & $iter$ & $time$ & $mv$ & $iter$ & $time$ & $mv$
      & $iter$ & $time$ & $mv$\\ \hline
      IRRHLB &114&6.81&811&63&6.56&769&40&6.54&693   \\ \hline
      IRRLB  &171&9.47&1210&69&7.18&841&42&6.76&727  \\ \hline
      IRHLB  &194&5.81&1371&77&4.34&937&45&3.94&778 \\ \hline
      IRLB   &202&5.50&1427&82&4.77&997&47&4.05&812  \\ \hline
      IRLBA  &170&1.61&1196&72&0.98&871&43&0.81&739 \\ \hline
      IRLANB &282&5.73&1706&99&3.89&1103&56&3.43&910 \\ \hline
    \end{tabular}
    \end{center}
\end{table}

\begin{table}[htp]
\caption {dw2048 for $k=1, 3, 5, 10, tol=1e-6$}\label{tabledw6}
\begin{center}
    \begin{tabular}{|c|c|c|c|c|c|c|c|c|c|} \hline
      \multicolumn{10}{|c|}{$k=1$}\\ \hline
      $m$&  \multicolumn{3} {|c|} {30}& \multicolumn{3} {|c|} {40}
      & \multicolumn{3} {|c|} {50}\\  \hline
      Algorithms & $iter$ & $time$ & $mv$ & $iter$ & $time$
      & $mv$& $iter$ & $time$ & $mv$\\ \hline
      IRRHLB &93&24.8&2242&77&38.1&2776&64&37.7&2948   \\ \hline
      IRRLB  &156&36.9&4060&52&23.3&1876&46&30.0&2120  \\ \hline
      IRHLB  &236&59.6&6140&128&57.2&4612&82&47.9&3776 \\ \hline
      IRLB   &266&67.5&6920&145&62.9&5224&93&57.5&4282  \\ \hline
      IRLBA  &276&9.71&7180&148&7.54&5332&94&6.64&4328 \\ \hline
      IRLANB &406&64.2&10155&219&70.0&7670&142&60.4&6395 \\ \hline
      \multicolumn{10}{|c|}{$k=3$}\\ \hline
      $m$&  \multicolumn{3} {|c|} {30}& \multicolumn{3} {|c|} {40}
      & \multicolumn{3} {|c|} {50}\\  \hline
      Algorithms & $iter$ & $time$ & $mv$ & $iter$ & $time$
      & $mv$& $iter$ & $time$ & $mv$\\ \hline
      IRRHLB &100&22.7&2406&81&34.8&2760&67&40.9&2954   \\ \hline
      IRRLB  &117&26.2&2814&72&31.5&2454&46&32.8&2030  \\ \hline
      IRHLB  &209&40.8&5022&110&41.8&3746&71&42.7&3130 \\ \hline
      IRLB   &230&44.0&5526&120&45.2&4086&78&44.9&3438  \\ \hline
      IRLBA  &238&7.84&5718&124&6.21&4222&80&5.55&3526 \\ \hline
      IRLANB &280&29.3&6447&142&42.7&4693&92&35.5&3963 \\ \hline
      \multicolumn{10}{|c|}{$k=5$}\\ \hline
      $m$&  \multicolumn{3} {|c|} {30}& \multicolumn{3} {|c|} {40}
      & \multicolumn{3} {|c|} {50}\\  \hline
      Algorithms & $iter$ & $time$ & $mv$ & $iter$ & $time$
      & $mv$& $iter$ & $time$ & $mv$\\ \hline
      IRRHLB &107&22.6&2362&80&30.9&2568&69&49.8&2906   \\ \hline
      IRRLB  &146&32.7&6220&75&35.2&2408&49&33.6&2066  \\ \hline
      IRHLB  &186&31.5&4100&95&32.7&3048&63&35.1&2654 \\ \hline
      IRLB   &266&47.0&5860&134&43.0&4296&85&47.8&3578  \\ \hline
      IRLBA  &287&9.27&6321&142&6.79&4552&89&6.06&3746 \\ \hline
      IRLANB &224&28.9&4713&114&23.4&3543&73&28.0&3002 \\ \hline
      \multicolumn{10}{|c|}{$k=10$}\\ \hline
      $m$&  \multicolumn{3} {|c|} {30}& \multicolumn{3} {|c|} {40}
      & \multicolumn{3} {|c|} {50}\\  \hline
      Algorithms & $iter$ & $time$ & $mv$ & $iter$ & $time$ & $mv$
      & $iter$ & $time$ & $mv$\\ \hline
      IRRHLB &190&38.9&3243&110&46.2&2983&89&67.7&3306   \\ \hline
      IRRLB  &246&51.4&4195&136&65.2&3685&68&54.1&2529  \\ \hline
      IRHLB  &443&63.5&7544&187&59.3&5062&110&56.8&4083 \\ \hline
      IRLB   &483&69.0&8224&203&63.6&5494&117&60.0&4342  \\ \hline
      IRLBA  &335&9.99&5636&159&7.21&4250&90&5.85&3322 \\ \hline
      IRLANB &546&54.7&8750&222&43.7&5786&126&44.2&4550 \\ \hline
    \end{tabular}
    \end{center}
\end{table}

\begin{table}[htp]
\caption {pde2961 for $k=1, 3, 5, 10, tol=1e-6$}\label{tablepde6}
\begin{center}
    \begin{tabular}{|c|c|c|c|c|c|c|c|c|c|} \hline
      \multicolumn{10}{|c|}{$k=1$}\\ \hline
      $m$&  \multicolumn{3} {|c|} {30}& \multicolumn{3} {|c|} {40}
      & \multicolumn{3} {|c|} {50}\\  \hline
      Algorithms & $iter$ & $time$ & $mv$ & $iter$ & $time$
      & $mv$& $iter$ & $time$ & $mv$\\ \hline
      IRRHLB &127&44.1&3306&92&41.3&3316&85&78.9&3914   \\ \hline
      IRRLB  &230&77.7&5984&147&80.3&5296&66&59.7&3040  \\ \hline
      IRHLB  &371&116&9650&193&106&6952&122&105&5616 \\ \hline
      IRLB   &425&136&11054&226&133&8140&142&122&6536  \\ \hline
      IRLBA  &463&21.6&12042&238&16.6&8572&148&14.5&6812\\ \hline
      IRLANB &490&100&12255&252&103&8825&157&99.0&7070 \\ \hline
      \multicolumn{10}{|c|}{$k=3$}\\ \hline
      $m$&  \multicolumn{3} {|c|} {30}& \multicolumn{3} {|c|} {40}
      & \multicolumn{3} {|c|} {50}\\  \hline
      Algorithms & $iter$ & $time$ & $mv$ & $iter$ & $time$
      & $mv$& $iter$ & $time$ & $mv$\\ \hline
      IRRHLB &149&44.7&3582&113&62.2&3848&101&83.1&4450   \\ \hline
      IRRLB  &254&83.6&6102&149&90.3&5072&85&63.3&3746  \\ \hline
      IRHLB  &475&137&11406&239&131&8132&146&113&6430 \\ \hline
      IRLB   &537&143&12894&272&134&9254&167&139&7354  \\ \hline
      IRLBA  &581&27.9&13950&284&19.3&9662&172&16.3&7574 \\ \hline
      IRLANB &686&133&15785&339&111&11194&204&97.3&8779 \\ \hline
      \multicolumn{10}{|c|}{$k=5$}\\ \hline
      $m$&  \multicolumn{3} {|c|} {30}& \multicolumn{3} {|c|} {40}
      & \multicolumn{3} {|c|} {50}\\  \hline
      Algorithms & $iter$ & $time$ & $mv$ & $iter$ & $time$
      & $mv$& $iter$ & $time$ & $mv$\\ \hline
      IRRHLB &161&52.0&3550&127&78.6&4082&100&97.6&4208   \\ \hline
      IRRLB  &239&65.6&5266&110&66.2&3528&103&99.2&4334  \\ \hline
      IRHLB  &518&122&11404&248&114&7944&149&111&6266 \\ \hline
      IRLB   &575&113&12658&278&126&8904&165&116&6938  \\ \hline
      IRLBA  &579&25.1&12745&273&18.2&8743&164&15.2&6895\\ \hline
      IRLANB &604&99.7&12693&284&84.6&8813&169&71.5&6938 \\ \hline
      \multicolumn{10}{|c|}{$k=10$}\\ \hline
      $m$&  \multicolumn{3} {|c|} {30}& \multicolumn{3} {|c|} {40}
      & \multicolumn{3} {|c|} {50}\\  \hline
      Algorithms & $iter$ & $time$ & $mv$ & $iter$ & $time$ & $mv$
      & $iter$ & $time$ & $mv$\\ \hline
      IRRHLB &290&83.7&4943&176&80.7&4765&139&145&5156   \\ \hline
      IRRLB  &571&145&9720&189&108&5116&144&132&5341  \\ \hline
      IRHLB  &925&212&15738&373&163&10084&205&150&7598 \\ \hline
      IRLB   &1004&193&17081&403&157&10894&223&140&8264  \\ \hline
      IRLBA  &601&23.3&10157&258&15.9&6948&145&12.8&5355 \\ \hline
      IRLANB &1302&138&20846&502&119&13066&269&106&9698 \\ \hline
    \end{tabular}
    \end{center}
\end{table}

\begin{table}[htp]
\caption {plat1919 for $k=1, 3, 5, 10,
tol=1e-6$}\label{tableplat6}
\begin{center}
    \begin{tabular}{|c|c|c|c|c|c|c|c|c|c|} \hline
      \multicolumn{10}{|c|}{$k=1$}\\ \hline
      $m$&  \multicolumn{3} {|c|} {15}& \multicolumn{3} {|c|} {20}
      & \multicolumn{3} {|c|} {25}\\  \hline
      Algorithms & $iter$ & $time$ & $mv$ & $iter$ & $time$
      & $mv$& $iter$ & $time$ & $mv$\\ \hline
      IRRHLB &146&8.20&1610&129&15.0&2068&74&12.0&1558   \\ \hline
      IRRLB  &351&22.7&3865&200&23.0&3204&105&19.0&2209  \\ \hline
      IRHLB  &605&37.5&6659&286&27.0&4580&163&28.3&3427 \\ \hline
      IRLB   &671&41.6&7385&313&32.9&5012&183&29.3&3847  \\ \hline
      IRLBA  &727&12.7&8001&337&8.08&5396&191&6.15&4015 \\ \hline
      IRLANB &943&33.2&9435&425&27.9&6380&234&29.2&4685 \\ \hline
      \multicolumn{10}{|c|}{$k=3$}\\ \hline
      $m$&  \multicolumn{3} {|c|} {15}& \multicolumn{3} {|c|} {20}
      & \multicolumn{3} {|c|} {25}\\  \hline
      Algorithms & $iter$ & $time$ & $mv$ & $iter$ & $time$
      & $mv$& $iter$ & $time$ & $mv$\\ \hline
      IRRHLB &214&10.9&1932&209&17.7&2890&153&21.4&2913   \\ \hline
      IRRLB  &605&30.5&5451&179&15.3&2512&123&16.4&2343  \\ \hline
      IRHLB  &923&37.2&8313&379&31.9&5312&214&28.4&4072 \\ \hline
      IRLB   &1089&48.7&9807&451&37.8&6320&253&28.8&4813  \\ \hline
      IRLBA  &1277&17.7&11499&511&10.5&7160&273&8.44&5193 \\ \hline
      IRLANB &1428&46.4&11431&545&32.7&7092&286&25.6&5155 \\ \hline
      \multicolumn{10}{|c|}{$k=5$}\\ \hline
      $m$& \multicolumn{3} {|c|} {15}& \multicolumn{3} {|c|} {20}
      & \multicolumn{3} {|c|} {25}\\  \hline
      Algorithms & $iter$ & $time$ & $mv$ & $iter$ & $time$
      & $mv$& $iter$ & $time$ & $mv$\\ \hline
      IRRHLB &385&15.3&2703&237&19.7&2582&171&22.2&2915   \\ \hline
      IRRLB  &1983&86.0&13889&454&40.8&5456&266&35.8&4530  \\ \hline
      IRHLB  &$n.c$&-&-&793&48.2&9524&387&43.7&6587 \\ \hline
      IRLB   &$n.c$&-&-&869&64.0&10436&432&46.7&7352 \\ \hline
      IRLBA  &$n.c$&-&-&722&13.5&8647&352&9.77&5983 \\ \hline
      IRLANB &$n.c$&-&-&1150&54.5&12659&524&45.9&8393 \\ \hline
      \multicolumn{10}{|c|}{$k=10$}\\ \hline
      $m$& \multicolumn{3} {|c|} {20}& \multicolumn{3} {|c|} {25}
      & \multicolumn{3} {|c|} {30}\\  \hline
      Algorithms & $iter$ & $time$ & $mv$ & $iter$ & $time$ & $mv$
      & $iter$ & $time$ & $mv$\\ \hline
      IRRHLB &685&45.1&4808&368&45.5&4429&395&94.9&6728   \\ \hline
      IRRLB  &$n.c$.&-&-&1283&165&15409&490&116&8343  \\ \hline
      IRHLB  &$n.c$.&-&-&1809&192&21721&891&155&15160 \\ \hline
      IRLB   &$n.c$.&-&-&$n.c$.&-&-&1428&209&24289  \\ \hline
      IRLBA  &$n.c$.&-&-&1317&26.8&12412&562&17.2&8330 \\ \hline
      IRLANB &$n.c$.&-&-&$n.c$.&-&-&1039&92.7&16638 \\ \hline
    \end{tabular}
    \end{center}
\end{table}

We found that all the algorithms computed the desired smallest
singular values correctly once they converged. The computed
smallest singular values had relative errors no more than a very
modest multiple of $\kappa(A)\times 10^{-6}$. We observed that, in
terms of restarts and matrix-vector products, IRRHLB was often
considerably more efficient and several times faster than the
others except IRRLB. IRRLB was nearly as efficient as IRRHLB in
many cases, and it was slightly better than IRRHLB in a few cases;
see, e.g., Table~\ref{tabledw6}--\ref{tablepde6} for the results
on dw2048 for $m=50$ and pde2961 for $k=1,3$, $m=50$. However, for
the relatively difficult plat1919, it was less robust than IRRHLB
and failed to converge for some $k$ and $m$; see
Table~\ref{tableplat6}. IRHLB, IRLB, IRLBA and IRLANB were less
robust and efficient than IRRLB, as the tables indicate. In
addition, the results demonstrate that the bigger $k$ was, the
more restarts the algorithms used generally.

For $opts.tol=1e-9$, we observed similar phenomena and had similar
findings. The only essential exception is that for plat1919 and
$k=5$, $m=15$, IRRLB did not converge after 2000 restarts were used.
Furthermore, we found that for the four test matrices all the
algorithms used more restarts for $opts.tol=1e-9$ than those for
$opts.tol=1e-6$ and they continued converging very smoothly from
$opts.tol=1e-6$ to $opts.tol=1e-9$, provided they converged. Hence
we do not list the results anymore.

\begin{figure}[t]
\epsfig{figure=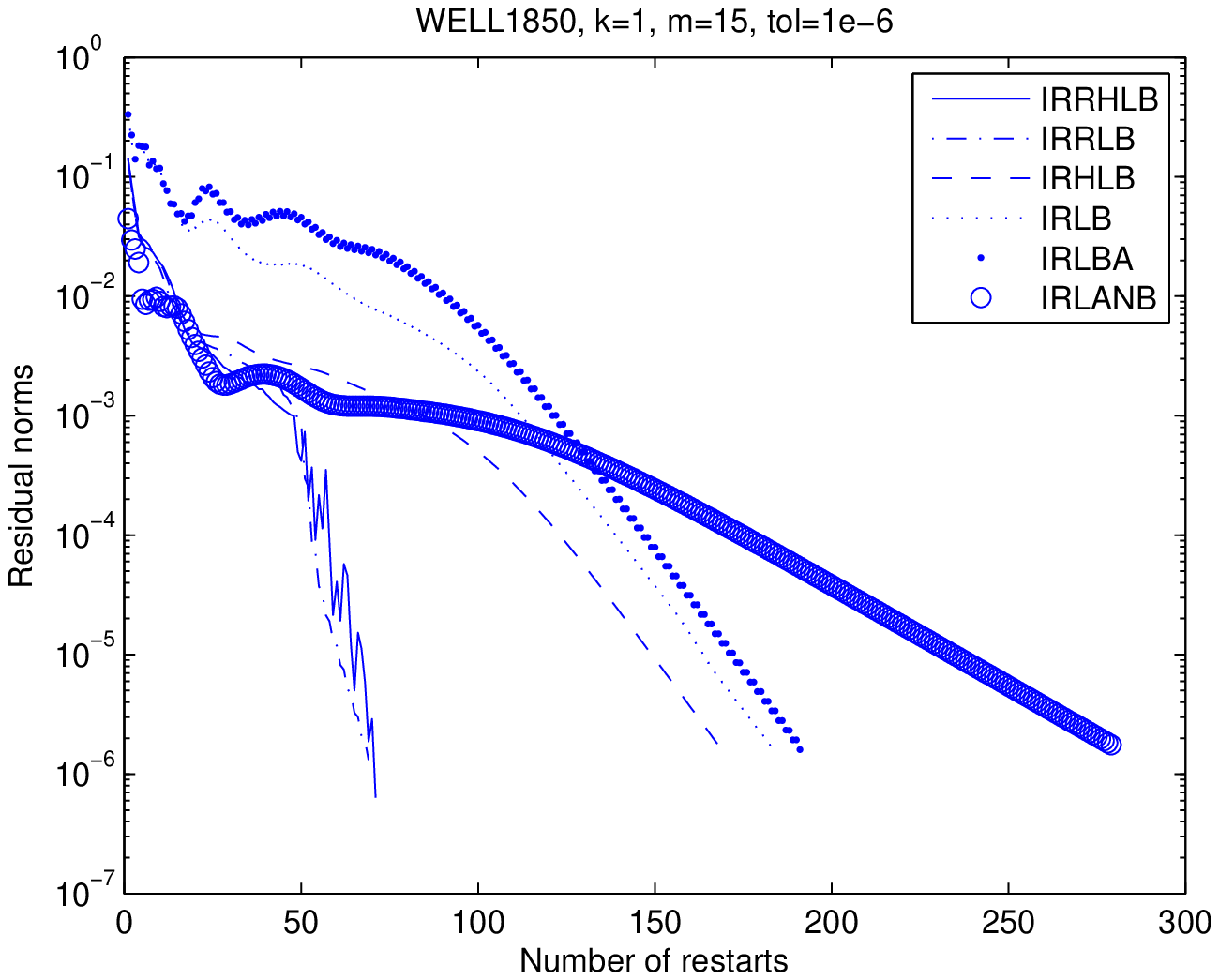,height=2.0in}
\epsfig{figure=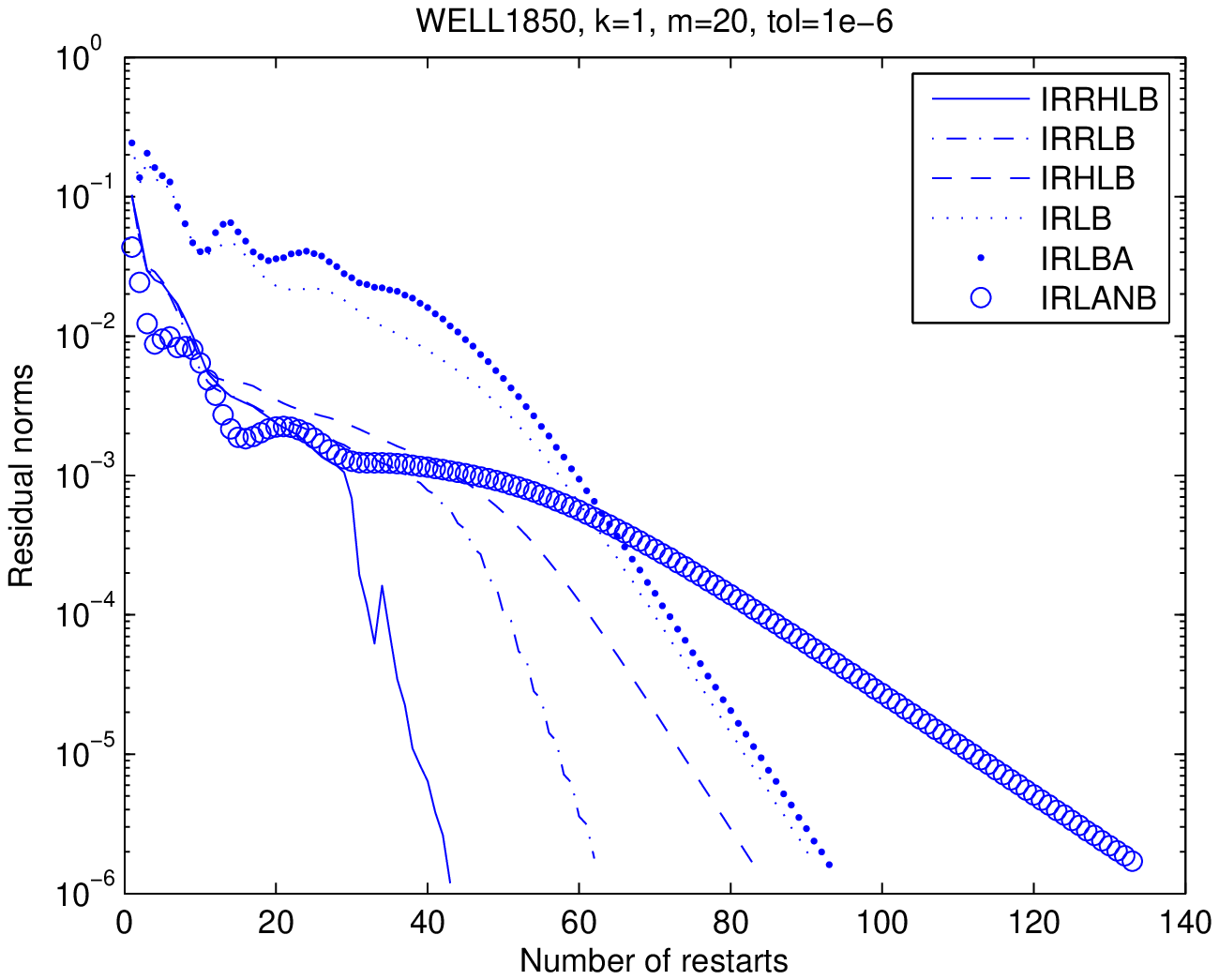,height=2.0in}
\epsfig{figure=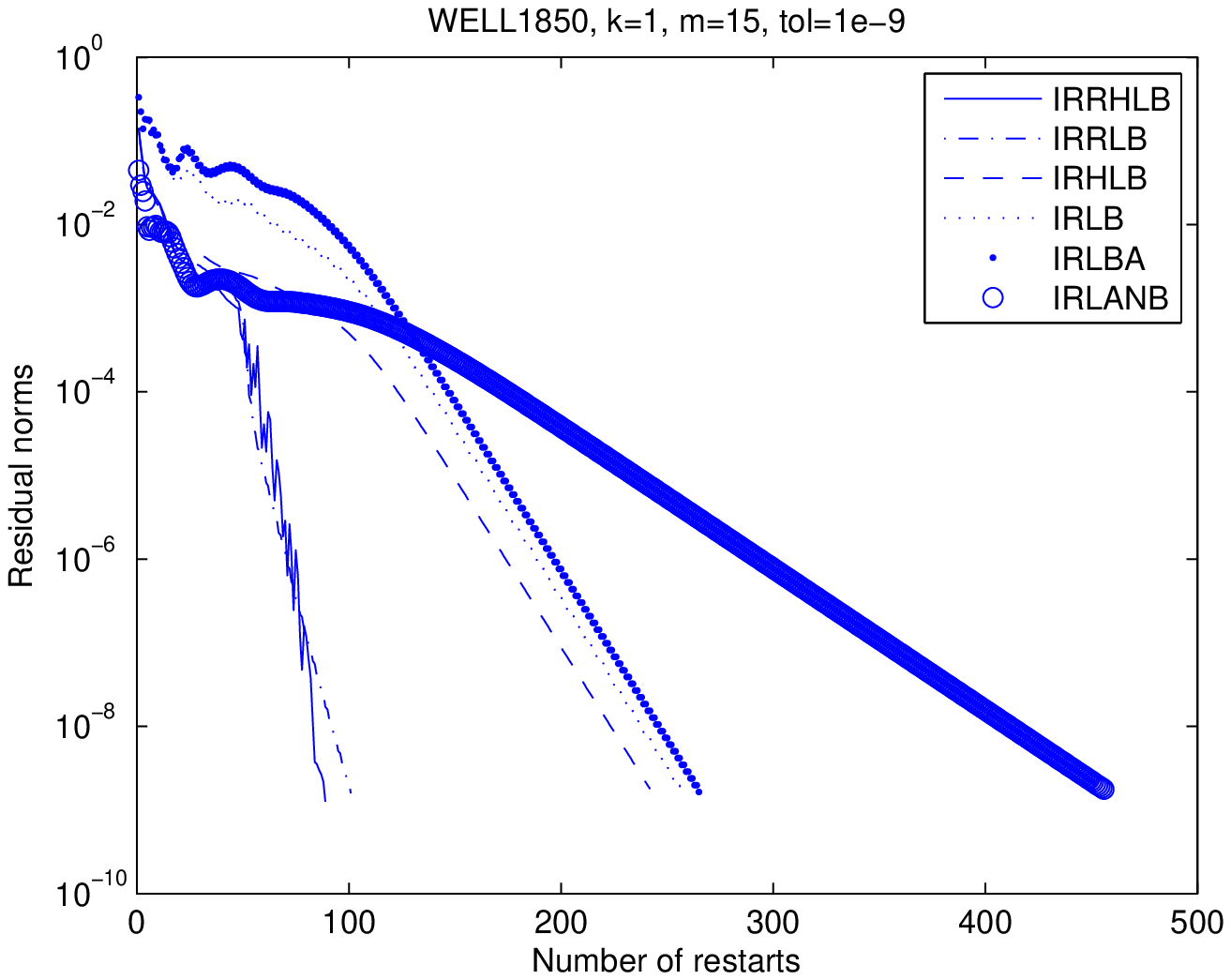,height=2.0in}
\epsfig{figure=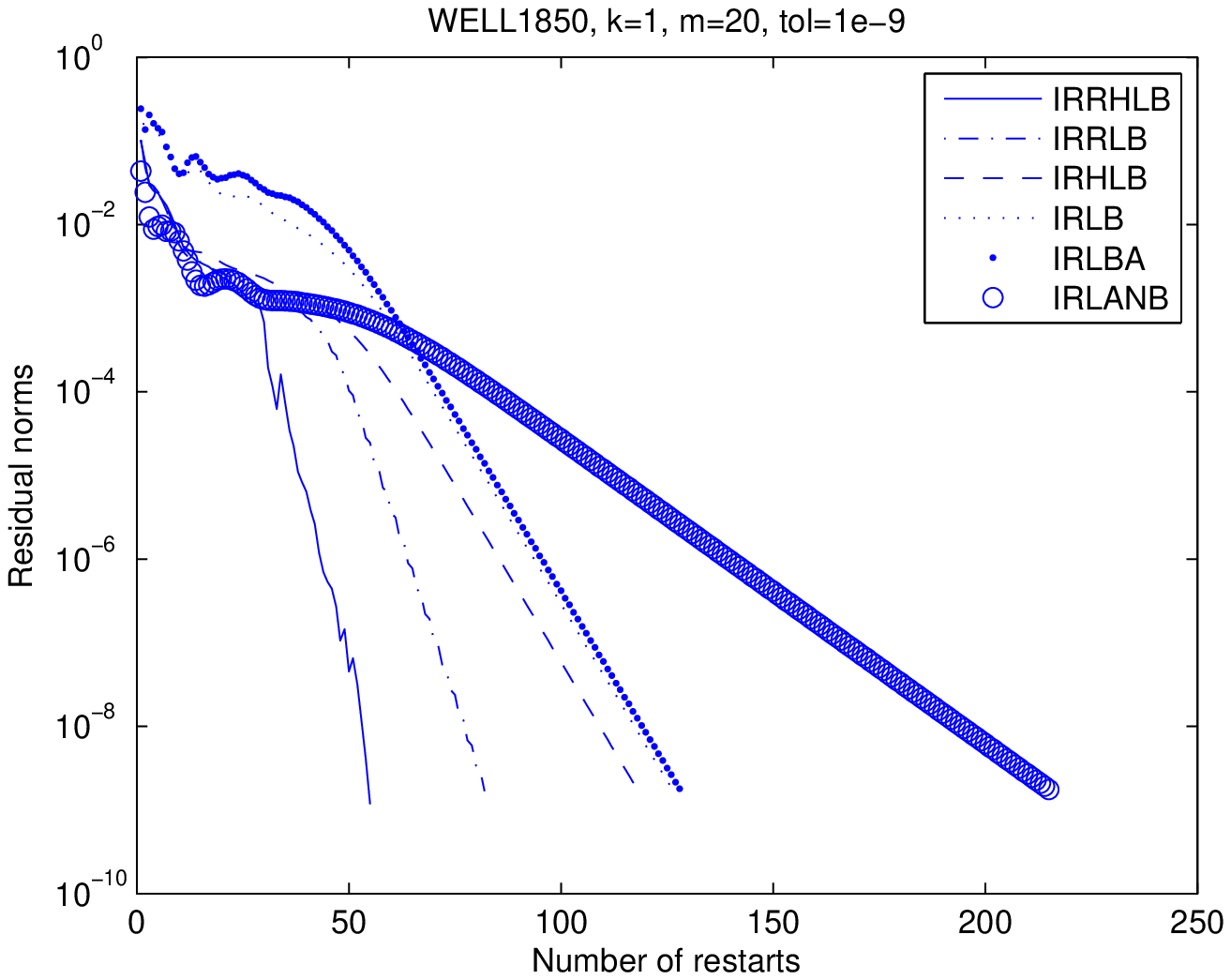,height=2.0in}
\caption{Convergence curves of well1850 with $k=1$ and
$tol=1e-6,\,1e-9$.}\label{fig3}
\end{figure}

\begin{figure}[t]
\epsfig{figure=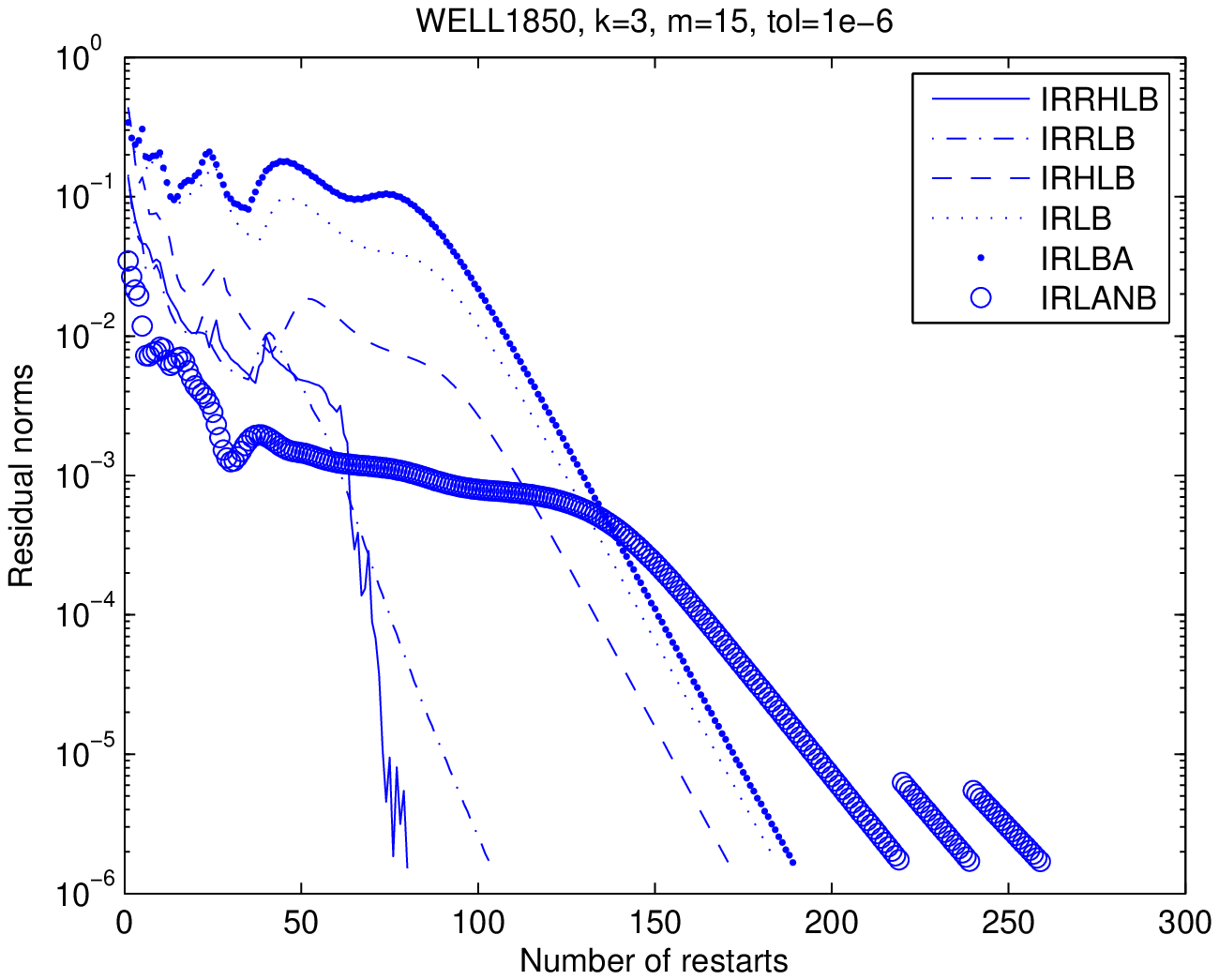,height=2.0in}
\epsfig{figure=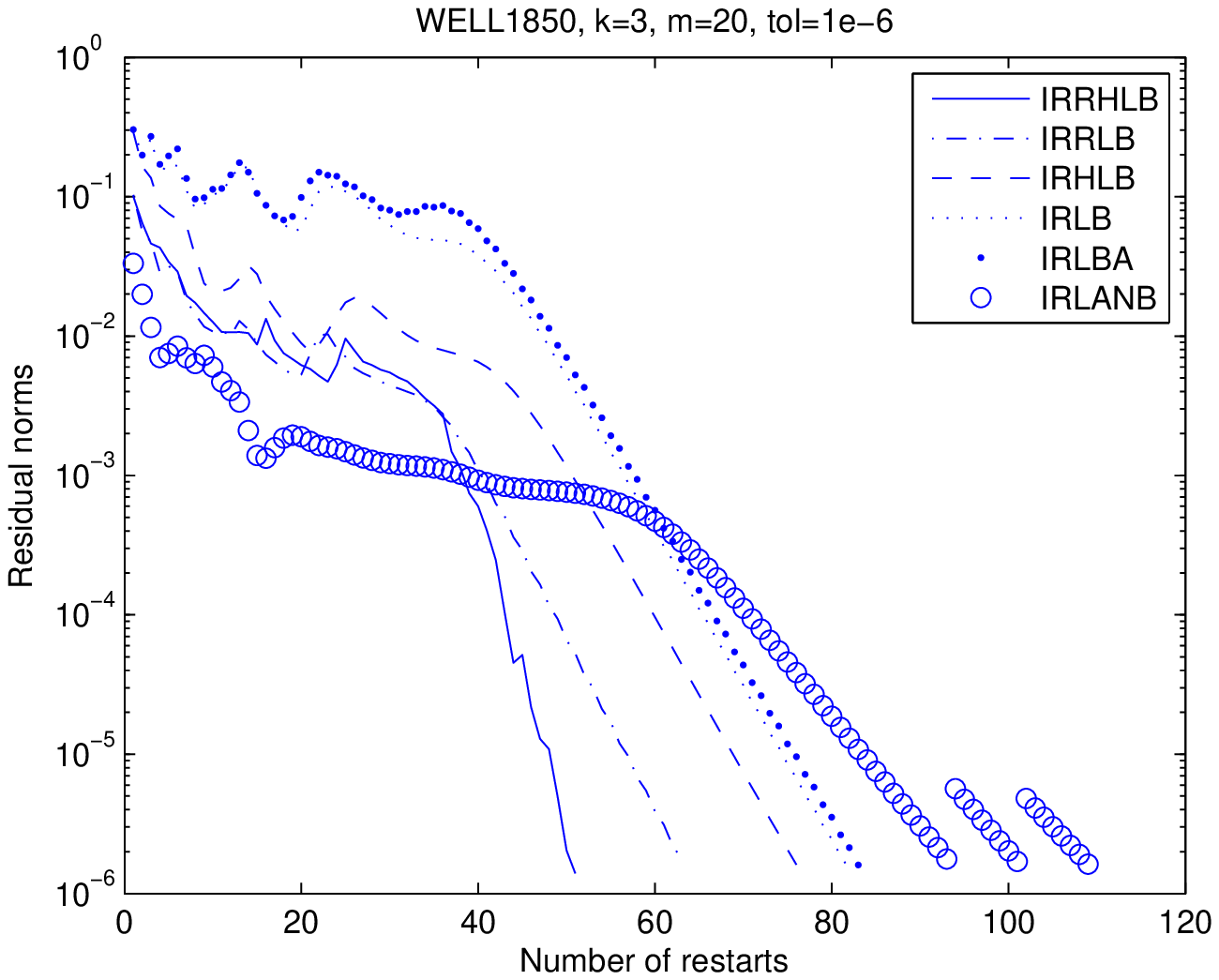,height=2.0in}
\epsfig{figure=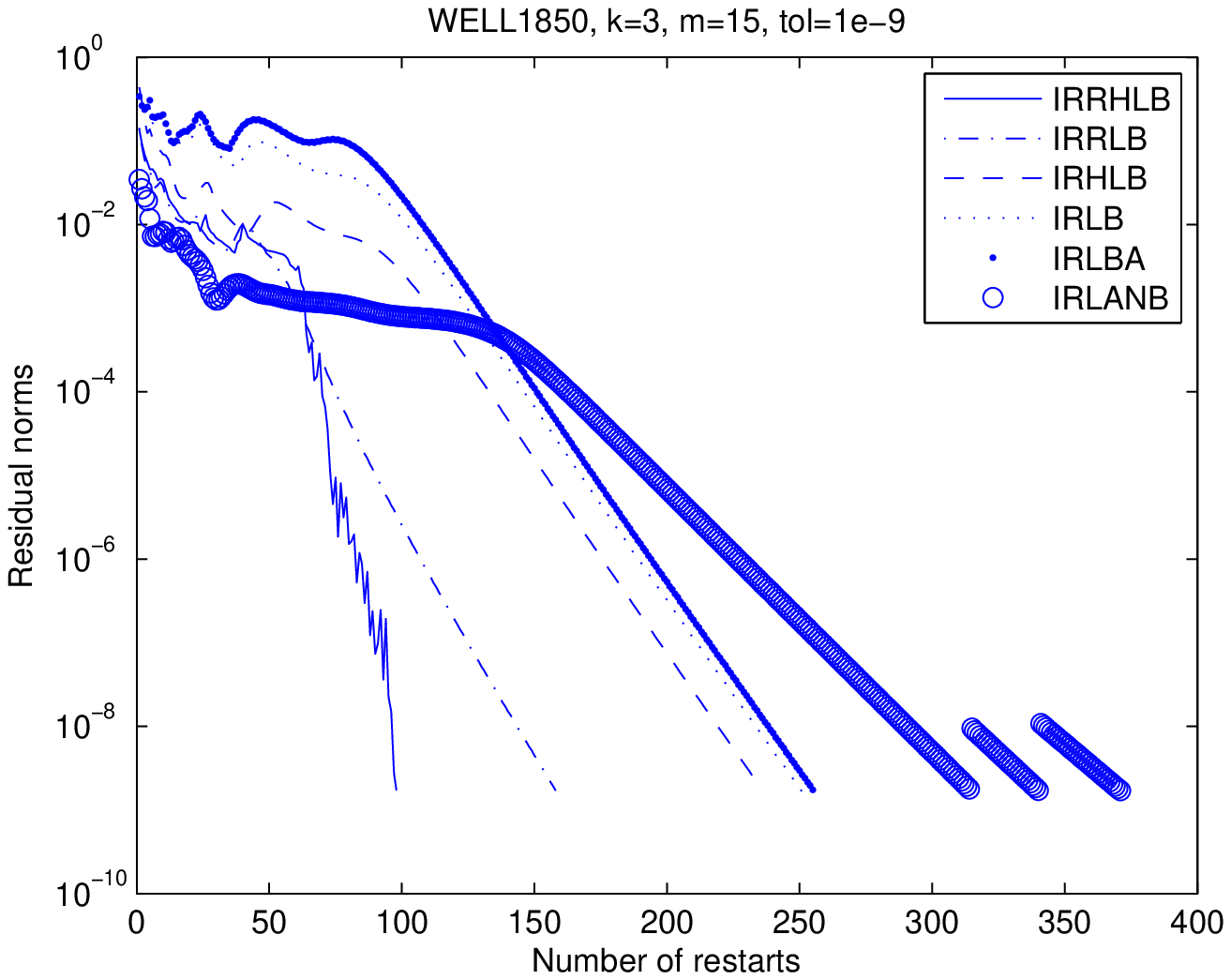,height=2.0in}
\epsfig{figure=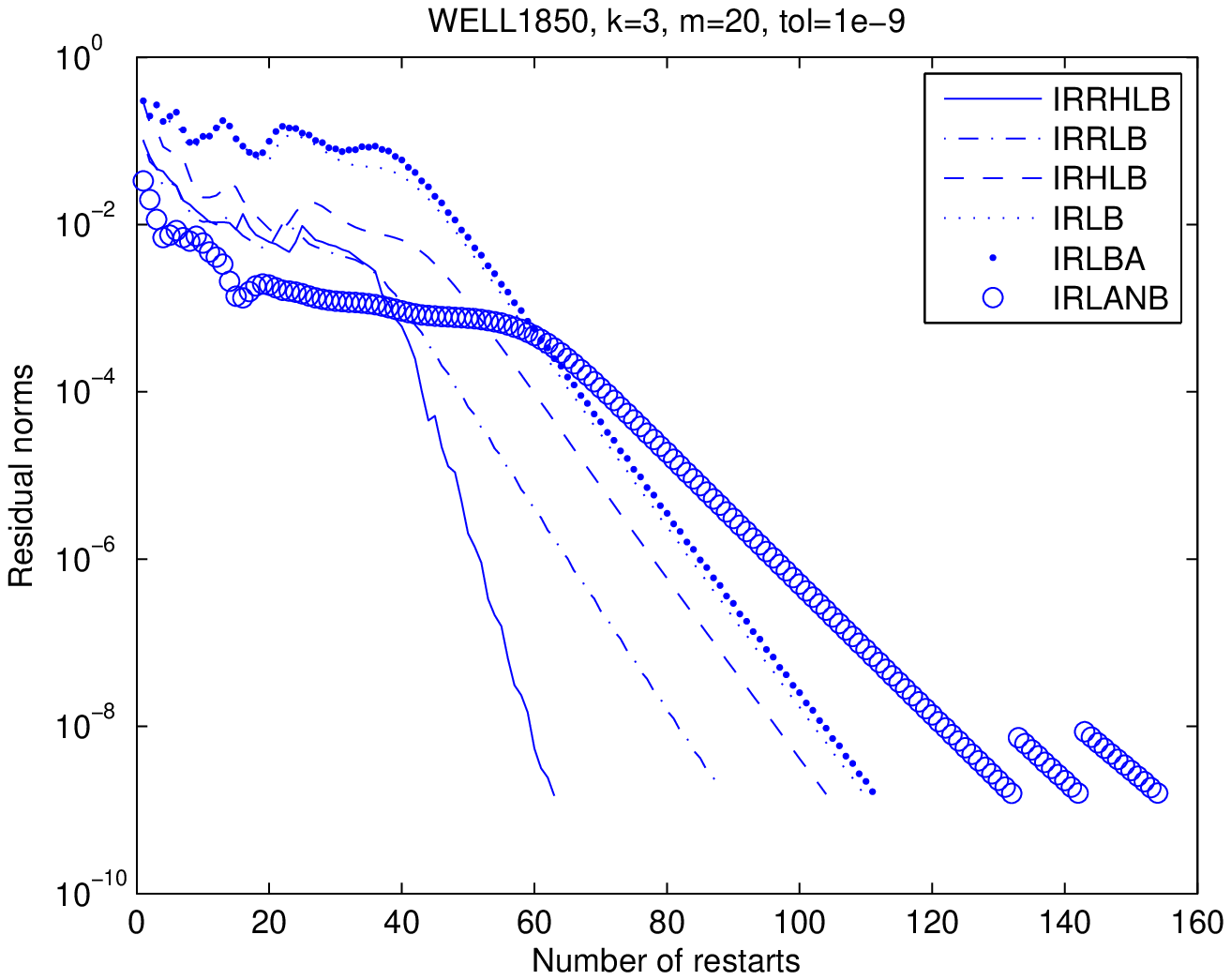,height=2.0in}
\caption{Convergence curves of well1850 with $k=3$ and $tol=1e-6,\
1e-9$.}\label{fig4}
\end{figure}

To be more illustrative, we draw some typical curves that feature
general convergence processes of the six algorithms.
Figures~\ref{fig3}--\ref{fig4} depict absolute residual norms
versus restarts for well1850 when $k=1,3$ and $opts.tol=1e-6,\
1e-9$, respectively. The figures clearly demonstrate that IRRHLB
is the fastest, IRRLB is the second best, IRHLB is faster than
IRLB while IRHLB, IRLB, IRLBA and IRLANB are comparable and
competitive though IRLANB may be slightly slower. The tables tell
us that IRHLB was faster than IRLB. We see from the figures that
after some stages the algorithms started converging quite smoothly
and they used more but not too more restarts for the smaller
$opts.tol=1e-9$. Besides, for IRLANB, we see that they computed
the smallest singular triplet after many restarts then found the
second and third smallest singular triplets very quickly. This is
because after the previous singular triplet(s) was (were) computed
the available subspaces had contained rich information on the
later desired singular vectors.

We have done more experiments and have similar findings. Based on
them, we may conclude that IRRHLB is the best and the most robust
for general purpose and IRRLB is the second best. As far as
overall efficiency is concerned, in terms of restarts and
matrix-vector products, IRRHLB is the fastest and IRRLB is the
second best while IRHLB, IRLB, IRLBA and IRLANB are all comparable
each other and no one is considerably superior to the others. A
further observation tells us that IRHLB is faster than IRLB. That
IRRHLB is superior to IRRLB and IRHLB is better than IRLB sheds
light on the fact argued in the introduction: The refined harmonic
projection and the harmonic projection are more suitable for
computing the smallest singular triplets than the refined standard
projection and the standard projection, respectively. Meanwhile,
we find that, as far as CPU timings are concerned, IRRHLB can be
inferior to IRLBA. This may be partly because $A$ is not very
large or $A$ too sparse, so that the savings of the first $k+3$
steps of the Lanczos bidiagonalization process cannot compensate
implicit restarting with $m-(k+3)$ shifts, and partly because our
code on implicit restarting is not far from optimized. In any
event, as a whole, we can draw an overall conclusion that IRRHLB
is at least competitive with and can be much more efficient than
the five other state of the art algorithms in both robustness and
efficiency.

The advantages of IRRHLB are twofold: It extracts the best left and
right approximate singular vectors from the given subspaces in the
sense of residual minimizations; it uses the better refined harmonic
shifts to construct better subspaces at each restart. Each of these
two advantages alone may not gain much, but, as restarts proceed,
the cumulative effect of their combination may be very striking, so
that IRRHLB can be much more efficient than the other algorithms, as
also noticed and commented on the refined algorithms for the large
eigenproblem in \cite{jia99a,jia02}.

\section{Concluding remarks}

We have presented the refined harmonic Lanczos bidiagonalization
method for computing the smallest singular triplets of large
matrices. We have developed a practical implicitly restarted
algorithm with the refined harmonic shifts scheme suggested. We
have done many numerical experiments and have compared the new
algorithm with the five other state of the art algorithms. The
results show that the new algorithm is at least competitive with
and can be much more efficient than the five other algorithms in
both robustness and efficiency.

We have reported the numerical results of computation of the
smallest singular triplets. We have also done many numerical
experiments on computation of the largest singular triplets. As
indicated in \cite{jianiu03}, IRRLB is generally preferable to
IRLB \cite{jianiu03,larsen}, PROPACK \cite{larsen1}, LANSO
\cite{larsen,larsen1} and {\sf svds} as well as some others; it is
the most robust among the restarted algorithms. Note that IRLANB
is designed to only compute the smallest singular triplets. For
computation of the largest singular triplets, we have found that
IRRLB is at least competitive with the four other algorithms, in
which IRLBA uses Ritz approximations. However, more observations
reveal that IRRHLB and IRHLB are considerably inferior to IRRLB
and IRLB, respectively. This suggests that IRRLB and IRLB are
suitable for computing both the largest singular triplets and the
smallest ones but IRRHLB and IRHLB are more suitable for computing
the smallest singular triplets.

The Matlab codes of IRRHLB, IRRLB, IRHLB and IRLB can be obtained
from the authors upon request.

\section*{Acknowledgements} We thank two referees very much for their
very valuable and helpful suggestions and comments, which made us
improve on the presentation considerably. Many thanks also go to
Kokiopoulou, Bekas, Gallopoulos and Baglama and Reichel for
generously providing us their IRLANB and IRLBA codes, which made our
numerical experiments and comparisons possible.

\end{document}